\documentclass[11pt]{article}
\usepackage{fullpage}
\usepackage{amsmath,amssymb,amsfonts,amsthm}
\usepackage{epsfig}
\usepackage{enumitem}
\usepackage{xcolor}
\usepackage{tikz, ifthen}
\usepackage[hidelinks]{hyperref}
\usepackage{thm-restate}
\usepackage{authblk}

\usetikzlibrary{quotes, positioning, arrows.meta}

\newtheorem{theorem}{Theorem}[section]
\newtheorem{lemma}[theorem]{Lemma}

\newtheorem{proposition}[theorem]{Proposition}

\newtheorem{corollary}[theorem]{Corollary}

\newtheorem{claim}{Claim}

\newenvironment{proofclaim}[1][{\it Proof of claim. \hspace{0.066cm}}]
{\noindent {}{#1}{}}{\QED\vspace{2ex}}

\def \QED {\hfill $\triangle$}      
\def \QD1 {\hfill $\spadesuit$}
\newcommand{\case}[2]{\smallskip {\bf Case #1\/:} {\it #2}}

\newcommand{\DF}[1]{{\bf #1\/}}

\newcommand{\etal}{et al.\,}

\newcommand{\set}[2]{\{#1 \;|\; #2 \}}
\newcommand{\ems}{\varnothing}
\newcommand{\sm}{\setminus}

\newcommand{\De}{\Delta}
\newcommand{\de}{\delta}

\newcommand{\el}{\ell}

\newcommand{\stil}{\scriptstyle}
\newcommand{\overr}[1]{\overset{\text{\tiny$\rightarrow$}}{#1}}
\newcommand{\dcn}[1]{\overr{\chi}(#1 : \cD)}
\newcommand{\ldcn}[1]{\overr{\chi}_{_{\stil\ell}}(#1: \cD)}
\newcommand{\pdcn}[1]{\overr{\chi}_{_{\rm DP}}(#1: \cD)}
\newcommand{\acn}{\overr{\chi}}
\newcommand{\lacn}{\overr{\chi}_{_{\stil\ell}}}
\newcommand{\pacn}{\overr{\chi}_{_{\rm DP}}}

\newcommand{\cn}{\chi}
\newcommand{\lcn}{\chi_{_{\stil\ell}}}
\newcommand{\pcn}{\chi_{_{\rm DP}}}

\newcommand{\cDG}{{\cal DG}}
\newcommand{\cD}{{\cal D}}
\newcommand{\cO}{{\cal O}}

\newcommand{\cAD}{{\cal AD}}
\newcommand{\cSD}{{\cal SD}}

\newcommand{\cBR}{{\cal BR}}

\newcommand{\cK}{\mathfrak{C}} 

\newcommand{\CR}{{\rm CR}}

\newcommand{\Crit}{\mathrm{Crit}}
\newcommand{\Crita}[1]{\mathrm{Crit(\cD, \acn,#1)}}
\newcommand{\Critb}[1]{\mathrm{Crit(\cD, \lacn,#1)}}
\newcommand{\Critc}[1]{\mathrm{Crit(\cD, \pacn,#1)}}

\newcommand{\Critaa}[1]{\mathrm{Crit(\cAD, \acn,#1)}}
\newcommand{\Critbb}[1]{\mathrm{Crit(\cAD, \lacn,#1)}}
\newcommand{\Critcc}[1]{\mathrm{Crit(\cAD, \pacn,#1)}}

\newcommand{\pa}{\xi}
\newcommand{\Critd}[1]{\mathrm{Crit(\cD, \pa,#1)}}

\newcommand{\dic}[1]{\vec{C}_{#1}}

\newcommand{\vf}{{\bf f}}

\newcommand{\f}{\varphi}
\newcommand{\fin}{\varphi^{-1}}

\newcommand{\dom}{{\rm dom}}
\newcommand{\su}{{\rm sp}}
\newcommand{\suc}{{\rm sp^o}}

\newcommand{\nat}{\mathbb{N}}
\newcommand{\nato}{\mathbb{N}_0}

\newcommand{\ganz}{\mathbb{Z}}

\newcommand{\NP}  {{\sf NP}}
\newcommand{\NPC} {{\NP}-complete}

\DeclareMathOperator{\UG}{UG}

\let\leq\leqslant
\let\geq\geqslant

\hypersetup{
	colorlinks,
	linkcolor={red!50!black},
	citecolor={blue!50!black},
	urlcolor={blue!80!black}
}

\newcommand{\jou}[4]{{\rm #1} #2 (#3) #4.}
\def \JCTB {J. Combin. Theory \, Ser.~B}
\def \DM {Discrete Math.}
\def \SIDMA {SIAM J. Discrete Math.}
\def \JGT {J. Graph Theory}

\numberwithin{equation}{section}

\begin{document}
\title{\bf Generalized DP-colorings of digraphs}

\author[,1]{Lucas Picasarri-Arrieta\thanks{Research supported by Japan Science and Technology Agency (JST) as part of Adopting Sustainable Partnerships for Innovative Research Ecosystem (ASPIRE), Grant Number JPMJAP2302. Email: {\tt lpicasarr@nii.ac.jp}}}
\author[,2]{Michael Stiebitz\thanks{Email: {\tt michael.stiebitz@tu-ilmenau.de}}}

\affil[1]{National Institute of Informatics, Tokyo, Japan}
\affil[2]{Technische Universit\"at Ilmenau, Ilmenau, Germany}

\date{}
\maketitle
\begin{abstract}
In this paper we consider the following three coloring concepts for digraphs. First of all, the generalized coloring concept, in which the same colored vertices of a digraph induce a subdigraph that satisfies a given digraph property. Second, the concept of variable degeneracy, introduced for graphs by Borodin, Kostochka and Toft in 2000; this allows to give a common generalization of the point partition number and the list dichromatic number. Finally, the DP-coloring concept as introduced for graphs by Dvořák and Postle in 2018, in which a list assignment of a graph is replaced by a cover. Combining these three coloring concepts leads to generalizations of several classical coloring results for graphs and digraphs, including the theorems of Brooks, of Gallai, of Erd\H{o}s, Rubin, and Taylor, and of Bernshteyn, Kostochka, and Pron for graphs, and the corresponding theorems for digraphs due to Harutyunyan and Mohar. Our main result combines the DP-coloring and variable degeneracy concepts for digraphs. 
\end{abstract}



\section{Introduction}

In this paper, we study the generalized coloring problem for digraphs, in which vertices of the same color satisfy a given digraph property. We also examine list and DP versions of this coloring concept.

The paper is organized as follows. Section~\ref{sect:preliminaries} contains some notation for graphs and digraphs, a brief overview of digraph properties, and the definition of three coloring parameters for a digraph $D$ with respect to a digraph property $\cD$: the $\cD$-dichromatic number $\dcn{D}$, the $\cD$-list-dichromatic number $\ldcn{D}$, and the $\cD$-DP-dichromatic number $\pdcn{D}$. In the study of coloring problems for digraphs, critical digraphs play an important role, since coloring problems for digraphs can often be reduced to coloring problems for critical digraphs. Section~\ref{sect:critical+main} presents some basic facts about critical digraphs and two of our main results of the paper, namely Theorems~\ref{theorem:mainDP} and \ref{theorem:mainLI}. We also discuss some interesting consequences of these results. A proof of Theorems~\ref{theorem:mainDP} and \ref{theorem:mainLI} is given in Section~\ref{sectsub:proofs:mainDP_LI}. 
The proof is based on Theorem~\ref{main_theorem}, another main result of the paper, which combines the DP-coloring concept with the concept of variable degeneracy for digraphs. We conclude the paper with some remarks in Section~\ref{sect:concluding}. 

\section{Preliminaries}
\label{sect:preliminaries}

Our notation is standard. In particular, $\nat$ denotes the set of positive integers and $\nato=\nat \cup \{0\}$. For integers $k$ and $\ell$, let $[k,\ell]=\set{x\in \ganz}{k \leq x \leq \ell}$. For two elements $p=(p_1,p_2)$ and $q=(q_1,q_2)$ of $\nato^2$, we denote by $p \leq q$ the relation ``$p_1 \leq q_1$ and $p_2 \leq q_2$". Additions on $\nato^2$ are made pairwise coordinate so $(p_1,q_1)+ (p_2,q_2) = (p_1 +p_2,q_1+q_2)$ and, for a nonnegative integer $k$, $(p_1,q_1) \cdot k = k \cdot (p_1,q_1) = (kp_1,kq_1)$.

\subsection{Graphs and digraphs}

We only consider finite simple graphs.
Let $G$ be a graph. Its order is denoted by $|G|$. If $|G|=0$, then $G$ is the \DF{empty graph}, and it is denoted $\ems$. If $G\not= \ems$, we denote by 
$$\de(G)=\min_{v\in V(G)}d_G(v) \mbox{~~~and~~~} \De(G)=\max_{v\in V(G)}d_G(v)$$ the \DF{minimum degree} and the \DF{maximum degree} of $G$, respectively. For $G=\ems$, we define $\de(G)=\De(G)=0$. 
 
 \smallskip
 
Our notation on directed graphs (digraphs for short) follow~\cite{Bang2009}. We only consider finite digraphs without any parallel arcs and loops.  Let $D$ be a digraph. Its \DF{order} is denoted by $|D|$. If $|D|=0$, then $D$ is the \DF{empty digraph}, and it is denoted $\ems$.
A \DF{matching} in $D$ is a set of pairwise disjoint arcs of $D$. A \DF{perfect matching} of $D$ is a matching of $D$ having $|D|/2$ arcs.
Let $X,Y\subseteq V(D)$, then $A_D(X,Y)$ denotes the set of arcs that have their initial vertex in $X$ and their terminal vertex in $Y$. So $X$ is \DF{independent} in $D$ if $A_D(X,X)=\ems$.

We call $C=(v_1,v_2,\ldots, v_n,v_1)$ a \DF{directed cycle} if $C$ is a directed graph with
$V(C)=\{v_1,v_2,\ldots,v_n\}$ and $A(C)=\{v_1v_2,v_2v_3, \ldots, v_{n-1}v_n, v_nv_1\}$ where the $v_i$ are pairwise distinct. 
The directed cycle of order $n$ is denoted by $\dic{n}$. A \DF{digon} is a pair of arcs with opposite directions between the same vertices. An arc of $D$ that is not in a digon of $D$ is called a \DF{simple arc} of $D$.

The digraph $D$ is a \DF{bidirected graph} if $D$ can be obtained from a graph $G$ by replacing each edge of $G$ by a digon, in this case we write $D=D^\pm(G)$.
The \DF{underlying graph} $\UG(D)$ of $D$ is the graph with vertex set $V(D)$ where $\{v,w\}$ is an edge of $\UG(D)$ if and only if at least one of $vw$ and $wv$ is an arc of $D$. 
We say that $D$ is \DF{connected} if $\UG(D)$ is connected. A \DF{component} of $D$ is a maximal subdigraph $D'$ of $D$ such that $D'$ is connected. 

If $D$ is connected, a \DF{separating vertex} of $D$ is a vertex $v \in V(D)$ such that $D-v$ has at least two components.
A digraph is \DF{2-connected} if it has at least three vertices and does not admit a separating vertex.
Furthermore, a \DF{block} of a digraph $D$ is a maximal connected subdigraph $D'$ of $D$ such that $D'$ has no separating vertex. Note that $D$ itself is a block if it is a connected digraph and has no separating vertex.

Let $a=vw$ be an arc of $D$. We say that $w$ is an \DF{out-neighbor} of $v$ in $D$ and that $v$ is an \DF{in-neighbor} of $w$ in $D$.
By $N_D^+(v)$ we denote the set of out-neighbors of $v$ in $D$; by $N_D^-(v)$ the set of in-neighbors of $v$ in $D$. The \DF{out-degree} $d^+_D(v)$ of $v$ is its number of out-neighbors. Similarly, its \DF{in-degree} $d^-_D(v)$ is its number of in-neighbors.
We say that $v$ is a \DF{source} in $D$ if $d_D^-(v)=0$, and it is a  \DF{sink} if $d_D^+(v)=0$. 
An \DF{antidirected cycle} is a digraph whose underlying graph is a cycle and in which each vertex is either a source or a sink.
The \DF{degree} of $v\in V(D)$ is the ordered pair $d_D(v) = (d^+_D(v),d^-_D(v))$. The digraph $D$ is \DF{Eulerian} if every vertex $u$ of $D$ satisfies $d^+_D(u) = d^-_D(u)$. It is further \DF{$r$-diregular} if  every vertex $u$ of $D$ satisfies $d^+_D(u) = d^-_D(u)=r$. If $D$ is nonempty, then let
\begin{itemize}
  \item $\De^+(D)=\max_{v\in V(D)} d_D^+(v)$ be the \DF{maximum out-degree} of $D$,
  \item $\de^+(D)=\min_{v\in V(D)} d_D^+(v)$ be the \DF{minimum out-degree} of $D$,
  \item $\De^-(D)=\max_{v\in V(D)} d_D^-(v)$ be the \DF{maximum in-degree} of $D$, and
  \item $\de^-(D)=\min_{v\in V(D)} d_D^-(v)$ be the \DF{minimum in-degree} of $D$.
\end{itemize}
By convention, we further set $\De^+(\ems)=\De^-(\ems)=\de^+(\ems)=\de^-(\ems)=0$.
We further define the \DF{maximum degree} $\De(D)$ and the \DF{minimum degree} $\de(D)$ of $D$ respectively as the ordered pairs
\[
	\De(D) = (\De^+(D), \De^-(D)) \mbox{~~~and~~~} \de(D) = (\de^+(D), \de^-(D)).
\]
 For a digraph $D$, we denote by $a_D:V(D)^2 \to \nato$ the \DF{adjacency function} of $D$, that is, for $u,v\in V(D)$ we have
$$a_D(u,v)=
\left\{ \begin{array}{ll}
1 & \mbox{\rm if } uv\in A(D),\\
0 & \mbox{\rm otherwise.}
\end{array}
\right.$$

\subsection{Digraph properties}

Let $\cDG$ denote the class of all digraphs. A \DF{digraph property} is a subclass of $\cDG$ that is closed with respect to isomorphisms. A digraph property $\cD$ is said to be \DF{nontrivial} if $\cD$ contains a nonempty digraph, but not all digraphs; otherwise $\cD$ is called \DF{trivial}. We call $\cD$ \DF{monotone} if $\cD$ is closed under taking subdigraphs, and we call $\cD$ \DF{hereditary} if $\cD$ is closed under taking induced subdigraphs. Clearly, every monotone digraph property is hereditary, but not conversely.
If $\cD$ is closed under taking (vertex) disjoint unions, then $\cD$ is called \DF{additive}. 

Let us consider some special digraph properties. A digraph is \DF{acyclic} if it contains no directed cycle. Let $\cAD$ denote the class of acyclic digraphs. It is easy to see that $\cAD$ is a digraph property that is nontrivial, monotone, and additive. Let $m\in \nato$ be a given integer. A digraph $D$ is called \DF{strictly $m$-degenerate} if every nonempty subdigraph $D'$ of $D$ has a vertex $v$ such that
$$\min\{d_{D'}^+(v),d_{D'}^-(v)\}<m.$$
Let $\cSD_m$ denote the class of digraphs that are strictly $m$-degenerate. Note that $\cSD_0=\{\ems\}$ is a trivial digraph property. Furthermore, it is well known and easy to prove that $\cAD=\cSD_1$. Clearly, for $m\geq 1$, the digraph property $\cSD_m$ is nontrivial, monotone, and additive. 

For a nontrivial hereditary digraph property $\cD$, let
$$\CR(\cD)=\set{D\in \cDG}{D \not\in \cD, \mbox{ but } D-v\in \cD \mbox{ for all } v\in V(D)}$$
and, moreover, define $d(\cD)=(d^+(\cD),d^-(\cD))$, where
$$d^+(\cD)=\min\set{\de^+(D)}{D\in \CR(\cD)} \mbox{~~~and~~~} d^-(\cD)=\min\set{\de^-(D)}{D\in \CR(\cD)}.$$
Note that $\CR(\cSD_m)$ contains the bidirected complete graph $D^\pm(K_{m+1})$, from which we easily conclude that $d(\cSD_m)=(m,m)$ for $m\geq 1$. We have the following general properties.

\begin{proposition}
\label{prop:smooth}
Let $\cD$ be a nontrivial and hereditary digraph property. Then the following statements hold:
\begin{enumerate}[label={\rm (\alph*)}]
\item $K_0, K_1 \in \cD$.
\label{enum:prop:smooth:a}
\item A digraph $D$ belongs to $\CR(\cD)$ if and only if each proper induced subdigraph of $D$ belongs to $\cD$, but $D$ itself does not belong to $\cD$.
\label{enum:prop:smooth:b}
\item $D\not\in \cD$ if and only if $D$ contains an induced subdigraph $D'$ with $D'\in \CR(\cD)$.
\label{enum:prop:smooth:c}
\item $\CR(\cD)\not= \ems$ and $d(\cD)\in \nato^2$.
\label{enum:prop:smooth:d}
\item If $D\not\in \cD$, but $D-v\in \cD$ for some vertex $v$ of $D$, then $d_D(v)\geq d(\cD)$.
\label{enum:prop:smooth:e}
\end{enumerate}
\end{proposition}

\begin{proof} Since $\cD$ is nontrivial, $\cD$ contains a nonempty digraph $D$. Since $\cD$ is hereditary, every induced subdigraph of $D$ belongs to $\cD$, which implies that both $K_0$ and $K_1$ belong to $\cD$. This proves~\ref{enum:prop:smooth:a}. Since $D-v$ is a proper induced subdigraph of a digraph $D$ for every vertex $v$ of $D$,~\ref{enum:prop:smooth:b} follows from the definition of $\CR(\cD)$ and the assumption that $\cD$ is hereditary.

To prove~\ref{enum:prop:smooth:c}, let $D$ be a digraph. If $D$ contains an induced subdigraph $D'$ with $D'\in \CR(\cD)$, then $D'\not\in \cD$ and thus $D\not\in \cD$ (because $\cD$ is hereditary). Now suppose that $D\not\in \cD$. Among all induced subdigraphs of $D$ not contained in $\cD$, let $D'$ be one whose order is minimum. Then $|D'|\geq 2$ (by~\ref{enum:prop:smooth:a}) and any proper induced subdigraph of $D'$ belongs to $\cD$. So it follows from~\ref{enum:prop:smooth:b} that $D'\in \CR(\cD)$. This completes the proof of~\ref{enum:prop:smooth:c}. Since $\cD$ is nontrivial, $\cD$ does not contain all digraphs. Hence,~\ref{enum:prop:smooth:c} implies~\ref{enum:prop:smooth:d}.

For the proof of~\ref{enum:prop:smooth:e}, suppose that $D\not\in \cD$, but $D-v\in \cD$ for some vertex $v$ of $D$. By~\ref{enum:prop:smooth:c}, $D$ contains an induced subdigraph $D'$ with $D'\in \CR(\cD)$. Then $v\in V(D')$, since otherwise $D'$ would be an induced subdigraph of $D-v$ and so $D'$ would belong to $\cD$ as $D-v$ belongs to $\cD$ and $\cD$ is hereditary, giving a contradiction. Then we deduce that
$d_D(v)\geq d_{D'}(v)\geq \de(D')\geq d(\cD)$. This proves~\ref{enum:prop:smooth:e}.
\end{proof}

\subsection{Generalized coloring of digraphs}

Let $D$ be a digraph. A \DF{coloring} of $D$ with \DF{color set} $C$ is a mapping $\f:V(D)\to C$. The sets $\fin(c)=\set{v\in V(D)}{\f(v)=c}$ with $c\in C$ are called \DF{color classes} of the coloring $\f$. A \DF{list assignment} of $D$ with \DF{color set} $C$ is a mapping $L:V(D) \to 2^C$ that assigns to each vertex $v\in V(D)$ a set (list) $L(v)\subseteq C$ of colors. If $|L(v)|\geq k$ for all $v\in V(D)$, we say that $L$ is a \DF{$k$-assignment} of $D$. A coloring $\f$ of $D$ is called an \DF{$L$-coloring} if $\f(v)\in L(v)$ for all $v\in V(D)$.

A \DF{cover} of $D$ is a pair $(X,H)$ consisting of a function $X$ and a digraph $H$ that is disjoint from $D$ satisfying the following two conditions:

\begin{enumerate}[label={(D\arabic*)}]
\item $X: V(D) \to 2^{V(H)}$ is a function that assigns to each vertex $v \in V(D)$ a vertex set $X_v = X(v) \subseteq V(H)$ such that the sets 
    $X_v$ with $v\in V(D)$ are pairwise disjoint.
\label{desc:D1}
\item $H$ is a digraph with vertex set $V(H)=\bigcup_{v\in V(D)}X_v$ such that each $X_v$ is an independent set of $H$. Furthermore, 
for any ordered pair $(u,v)$ of distinct vertices of $D$, $A_H(X_u,X_v)$ is a matching of $H$ (possibly empty) if $uv\in A(D)$ and $A_H(X_u,X_v)=\ems$ if $uv\not\in A(D)$. 
\label{desc:D2}
\end{enumerate}

The following useful result is an immediate consequence of the definition of a cover of a digraph.

\begin{proposition}
\label{prop:cover}
Let $D$ be a digraph, $(X,H)$ be a cover of $D$, $u$ and $v$ be distinct vertices of $D$, and $x_v\in X_v$. Then we have
$$\sum_{x\in X_u} a_H(x,x_v) \leq a_D(u,v) \mbox{~~~and~~~}  \sum_{x\in X_u} a_H(x_v,x) \leq a_D(v,u).$$
\end{proposition}

Let $(X,H)$ be a cover of $D$. If $|X_v|\geq k$ for all $v\in V(D)$, we say that $(X,H)$ is a \DF{$k$-cover} of $D$. A \DF{transversal} of $(X,H)$ is a vertex set $T\subseteq V(H)$ such that $|T \cap X_v|=1$ for all $v \in V(D)$. A set $T\subseteq V(H)$ is called \DF{partial transversal} of $(X,H)$ if $|T \cap X_v|\leq 1$ for all $v\in V(D)$. For $Y\subseteq V(H)$, let $\dom(Y:D)=\set{v\in V(D)}{X_v\cap Y\not=\ems}$ be the \DF{domain} of $Y$ in $D$.

Colorings of digraphs are only of interest if some restrictions are placed on the color classes. A natural restriction is to require that the subdigraphs induced by the color classes satisfy a given digraph property. However, such colorings make sense only for certain digraph properties.
A digraph property is called \DF{reliable} if it is nontrivial, hereditary and additive. A reliable digraph property that is also monotone is called \DF{strongly reliable}. Note that $\cSD_m$ with $m\geq 1$ is a strongly reliable digraph property. 

For the rest of this section, assume that $\cD$ is a reliable digraph property and $D$ is an arbitrary digraph. A \DF{$\cD$-coloring} of $D$ with color set $C$ is a coloring $\f$ of $D$ with color set $C$  such that $D[\fin(c)]\in \cD$ for all $c\in C$. If $L$ is a list assignment for $D$ with color set $C$,
then a \DF{$(\cD,L)$-coloring} of $D$ is an $L$-coloring $\f$ of $D$ such that
$D[\fin(c)]\in \cD$ for all $c\in C$. The \DF{$\cD$-dichromatic number} of $D$, denoted by $\dcn{D}$, is the least integer $k$ for which $D$ admits a $\cD$-coloring with a set of $k$ colors. The \DF{$\cD$-list-dichromatic number} of $D$, denoted by $\ldcn{D}$, is the least integer $k$ such that $D$ has a $(\cD,L)$-coloring whenever $L$ is a $k$-assignment of $D$. If $(X,H)$ is a cover of $D$, then a \DF{$\cD$-transversal} of $(X,H)$ is a transversal $T$ of $(X,H)$ such that $H[T]\in \cD$. A $\cD$-transversal of $(X,H)$ is also called a \DF{$(\cD,(X,H))$-coloring} of $D$. Taking $L(v)=X_v$ for all $v\in V(D)$, note that $D$ admits a $(\cD,(X,H))$-coloring if and only if $D$ has an $L$-coloring $\f$  such that $T=\set{\f(v)}{v \in V(D)}$ is a $\cD$-transversal of $(X,H)$. The \DF{$\cD$-DP-dichromatic number} of $D$, denoted by $\pdcn{D}$, is the least integer $k$ such that $D$ admits a $(\cD,(X,H))$-coloring whenever $(X,H)$ is a $k$-cover of $D$. We have the following inequalities for every digraph $D$, which in particular imply that the parameters above are well defined:
\begin{align} \label{equation:dcn<ldcn<pdcn}
\dcn{D}\leq \ldcn{D}\leq \pdcn{D} \leq 2|D|.
\end{align}
 The first inequality follows from the fact that a $\cD$-coloring of $D$ with color set $C$ may be considered as a $(\cD,L)$-coloring of $D$ for the constant list assignment $L\colon v\mapsto C$. 
 To see the second inequality, suppose that $\pdcn{D}=k$ and let $L$ be a list assignment for $D$ with $|L(v)|\geq k$ for all $v\in V(D)$. Define $(X,H)$ to be the cover of $D$ such that $X_v=\{v\}\times L(v)$ for all $v\in V(D)$ and, for two distinct vertices $(u,c)$ and $(v,c')$ of $H$, we have $(u,c)(v,c')\in A(H)$ if and only if $uv\in A(D)$ and $c=c'$.
We say that $(X,H)$ is the cover \DF{associated with the list assignment} $L$; in this case we write $(X,H)=C(D,L)$. It is easy to check that $(X,H)$ is indeed a $k$-cover of $D$. Furthermore, since $\cD$ is reliable, we conclude that $(X,H)$ has a $\cD$-transversal if and only if $D$ admits a $(\cD,L)$-coloring. This implies, in particular, that $\ldcn{D}\leq k$.
To see the last inequality, let $(X,H)$ be an arbitrary $k$-cover of $D$ with $k=2|D|$. Using induction on $|D|$, it follows from~\ref{desc:D2} that $(X,H)$ has a transversal $T$ such that $T$ is an independent set of $H$. We have $K_1\in \cD$ by Proposition~\ref{prop:smooth}\ref{enum:prop:smooth:a} and so $H[T]\in \cD$ (as $\cD$ is additive). Hence, $T$ is a $\cD$-transversal of $(X,H)$ implying that $\pdcn{D}\leq k=2|D|$. Thus \eqref{equation:dcn<ldcn<pdcn} is proved.

\smallskip

We need a few specific notation on covers of digraphs. Let $D$ be a nonempty digraph, and let $(X,H)$ be a cover of $D$. Then for every arc $uv\in A(D)$, the arc set $A_H(X_u,X_v)$ is a matching of $H$ and $A(H)$ is the union of all these matchings (by~\ref{desc:D2}). We say that $(X,H)$ is a \DF{saturated cover} of $D$ if $A_H(X_u,X_v)$ is a perfect matching of $H[X_u \cup X_v]$ for every arc $uv\in A(D)$.
For some integer $r\in \nat$, we further say that $(X,H)$ is \DF{$r$-uniform} if $|X_v|=r$ for all $v\in V(D)$. Observe that, if $(X,H)$ is saturated and $D$ is connected, then $(X,H)$ is $r$-uniform for some integer $r\in \nat$.

Let $U\subseteq V(H)$ be an arbitrary set. Let $D'=D[\dom(U:D)]$, let $H'=H[U]$, and let $X':V(H') \to 2^U$ be the map with $X'_v=X_v\cap U$ for all $v\in V(H')$. Then $(X',H')$ is a cover of $D'$ and it is denoted by $(X',H')=(X,H)/U$; we call this cover a \DF{subcover} of $(X,H)$ \DF{restricted} to $U$. Conversely, for any induced subdigraph $D'$ of $D$, we denote by $(X,H)/D'$ the subcover of $(X,H)$ restricted to $U=\bigcup_{v\in V(D')}X_u$.

Observe that, for every $k$-cover $(X,H)$ of $D$, there is a set $U\subseteq V(H)$ such that $(X',H')=(X,H)/U$ is a $k$-uniform cover of $D$. Clearly, if $T$ is a $\cD$-transversal of $(X',H')$, then $T$ is also a $\cD$-transversal of $(X,H)$. Hence, to show that $\pdcn{D}\leq k$ it suffices to show that every $k$-uniform cover $(X,H)$ of $D$ has a $\cD$-transversal. We will use this property many times along the paper.

We conclude this subsection with some simple but useful facts about the coloring parameters $\acn, \lacn$ and $\pacn$. 

\begin{proposition}
Let $\cD$ be a reliable digraph property and let $D$ be an arbitrary digraph. Then the following statements hold.
\begin{enumerate}[label={\rm(\alph*)}]
\label{prop:basicA}
\item If $\pa\in \{\acn, \lacn, \pacn\}$, then $\pa(D:\cD)=0$ if and only if $D=\ems$.
\label{enum:prop:chi:0and1:a}
\item If $\pa\in \{\acn, \lacn\}$, then $\pa(D:\cD)\leq 1$ if and only if $D\in \cD$.
\label{enum:prop:chi:0and1:b}
\item If $\cD$ is strongly reliable, then $\pdcn{D}\leq 1$ if and only if $D\in \cD$.
\label{enum:prop:chi:0and1:c}
\end{enumerate}
\end{proposition}
\begin{proof}
The proofs of~\ref{enum:prop:chi:0and1:a} and~\ref{enum:prop:chi:0and1:b} are straightforward. For the proof of~\ref{enum:prop:chi:0and1:c}, assume that $\cD$ is strongly reliable. We may assume that $D\not=\ems$. If $(X,H)$ is a 1-uniform cover of $D$, then $T=V(H)$ is the only transversal of $(X,H)$ and $H[T]=H\subseteq D$, and for $H$ we can choose any subdigraph of $D$ (by ~\ref{desc:D2}). Since $\cD$ is monotone, this implies that $\pdcn{D}\leq 1$ if and only if $D\in \cD$. 
\end{proof} 

Note that the assumption in statement~\ref{enum:prop:chi:0and1:c} is necessary. To see this, take for instance $\cD$ to be the class of digraphs that are disjoint union of bidirected complete graphs. Then $\cD$ is a reliable digraph property, but $\cD$ is not monotone, and $D=D^\pm(K_3)$ has a $1$-uniform cover $(X,H)$ such that $H$ is a $\dic{3}$ and does not belong to $\cD$. This shows that $\pdcn{D}\geq 2$ while $D\in \cD$. 

\begin{proposition}
\label{prop:basicB}
Let $\cD$ be a reliable digraph property, let $D$ be an arbitrary digraph, and let $\pa\in \{\acn, \lacn, \pacn\}$. 
Then the following statements hold.
\begin{enumerate}[label={\rm(\alph*)}]
\item If $D'$ is an induced subdigraph of  $D$, then $\pa(D':\cD) \leq \pa(D:\cD)$.
	\label{prop:basic_properties_reliable:a}
\item If $D\not=\ems$, then $\pa(D:\cD)=\max \set{\pa(D':\cD)}{D' \mbox{ is a component of } D}$.
\label{prop:basic_properties_reliable:b}
\item If $v\in V(D)$, then $\pa(D:\cD) \leq \pa(D-v:\cD) + 1$, unless $\pa=\pacn$ and $d(\cD)=(0,0)$.
\label{prop:basic_properties_reliable:c}
%
%
\item $\pa(D:\cD)\leq |D|$, unless $\pa=\pacn$ and $d(\cD)=(0,0)$. 
\label{prop:basic_properties_reliable:e}
\end{enumerate}
\end{proposition}
\begin{proof}
	Let us first prove~\ref{prop:basic_properties_reliable:a}. When $\pa \in \{\acn,\lacn\}$, this is a direct consequence from the fact that $\cD$ is hereditary. Let $(X',H')$ be an arbitrary $k$-cover of $D'$, where $k=\pacn(D:\cD)$. Observe that there is a $k$-cover $(X,H)$ of $D$ such that $(X',H')=(X,H)/D'$. Hence, there is a $\cD$-transversal $T$ of $(X,H)$. Then $T'=T\cap V(H')$ is a transversal of $(X',H')$ and $H'[T']=H[T']$ is an induced subdigraph of $H[T]$. Since $\cD$ is hereditary, this implies that $T'$ is a $\cD$-transversal of $(X',H')$, and so $\pdcn{D'}\leq k=\pdcn{D}$. This proves~\ref{prop:basic_properties_reliable:a}.
Statement~\ref{prop:basic_properties_reliable:b} follows easily from~\ref{prop:basic_properties_reliable:a} and the fact that $\cD$ is additive.

We now prove~\ref{prop:basic_properties_reliable:c}. We prove it formally for $\pa = \pacn$, the case $\pa  \in \{\acn,\lacn\}$ being similar. Assume that $d(\cD)\neq (0,0)$, hence $d^+(\cD)\geq 1$ or $d^-(\cD)\geq 1$. By directional duality, let us assume $d^+(\cD)\geq 1$. Let $k=\pdcn{D-v}$ and let $(X,H)$ be an arbitrary $(k+1)$-cover of $D$.
Let $x\in X_v$ and let $(X',H')$ be the cover of $D'=D-v$ such that $X_u'=X_u\sm N^+_H(x)$ for all $u\in V(D')$ and $H'=H-(X_v \cup N^+_H(x))$. 
By \ref{desc:D2}, $(X',H')$ is a $k$-cover of $D'$ and, therefore, $(X',H')$ has a $\cD$-transversal $T'$, that is, $T'$ is a transversal of $(X',H')$ and $H'[T']\in \cD$. Then $T=T'\cup \{x\}$ is a transversal of $(X,H)$.
If $H[T]\not \in \cD$, then we have $H[T]-x=H[T']=H'[T']\in \cD$ and so 
$d^+_{H[T]}(x)\geq d^+(\cD)\geq 1$ (by Proposition~\ref{prop:smooth}\ref{enum:prop:smooth:e}), a contradiction. Hence $(X,H)$ has $\cD$-transversal implying that
$\pdcn{D}\leq k+1=\pdcn{D-v}+1$. This proves~\ref{prop:basic_properties_reliable:c}. 

Statement~\ref{prop:basic_properties_reliable:e} follows from~\ref{prop:basic_properties_reliable:c} and Proposition~\ref{prop:basicA}\ref{enum:prop:chi:0and1:a} by induction on $|D|$. 
\end{proof} 

\subsection{\texorpdfstring{$\cAD$}{AD}-coloring of digraphs versus proper coloring of graphs}

Generalized colorings of graphs, where the same colored vertices induce a subgraph that satisfies a given graph property, are a well-established part of graph coloring theory; the reader will find a long list of relevant references in \cite[Sections 2.6 and 2.8]{StiebitzST2024}. The best studied graph property in the context of generalized colorings is the class $\cO$ of edgeless graphs. Let $G$ be an arbitrary graph. An \DF{$\cO$-coloring} of $G$ with \DF{color set} $C$ is a map $\f:V(G) \to C$ such that $G[\fin(c)]\in \cO$ for every color $c\in C$; the usual name for $\cO$-coloring is \DF{proper coloring}. The \DF{chromatic number} of $G$, denoted by $\cn(G)$, is the least integer $k$ for which $G$ admits an $\cO$-coloring with a set of $k$ colors. If $L$ is a list assignment of $G$ with color set $C$, then an \DF{$(\cO,L)$-coloring} of $G$ with \DF{color set} $C$ is an $\cO$-coloring $\f$ of $G$ with color set $C$ such that $\f(v)\in L(v)$ for every vertex $v\in V(G)$. The \DF{list-chromatic number} of $G$, denoted by $\lcn(G)$ is the least integer $k$ such that $G$ has an $(\cO,L)$-coloring for every $k$-assignment $L$ of $G$ (i.e., $L$ is a list assignment of $G$ with $|L(v)|\geq k$ for every vertex $v\in V(G)$). A \DF{cover} of $G$ is a pair $(X,H)$ consisting of a function $X$ and a graph $H$ (disjoint from $G$) such that the following conditions hold: 

\begin{enumerate}[label=(G\arabic*)]
\item $X: V(G) \to 2^{V(H)}$ is a function that assigns to each vertex $v \in V(G)$ a vertex set $X_v = X(v) \subseteq V(H)$ such that the sets $X_v$ with $v\in V(D)$ are pairwise disjoint.
\item $H$ is a graph with vertex set $V(H) = \bigcup_{v \in V(G)}X_v$ such that each $X_v$ is an independent set of $H$. Furthermore, for any pair $(u,v)$ of distinct vertices of $G$, $E_H(X_u,X_v)$ is a matching of $H$ (possibly empty) if $uv\in E(G)$ and $E_H(X_u,X_v)=\ems$ if $uv\not\in E(G)$. Here $E_H(X,Y)=\set{uv\in E(G)}{u\in X, v\in Y}$. 
\end{enumerate}

Let $(X,H)$ be a cover of $G$. We say that $(X,H)$ is a \DF{$k$-cover} of $G$ if $|X_v|\geq k$ for every $v\in V(G)$.  A \DF{transversal} of $(X,H)$ is a vertex set $T\subseteq V(H)$ such that $|T\cap X_v|=1$ for every $v\in V(G)$. An \DF{$\cO$-transversal} of $(X,H)$ is a transversal $T$ of $(X,H)$ such that $H[T]\in \cO$. An $\cO$-transversal of $(X,H)$ is also called an \DF{$(\cO,(X,H))$-coloring} of $G$. The \DF{DP-chromatic number} of $G$, denoted by $\pcn(G)$, is the least integer $k$ such that $G$ admits an $(\cO,(X,H))$-coloring whenever $(X,H)$ is a $k$-cover of $G$. Analogously to~\eqref{equation:dcn<ldcn<pdcn}, we have
\begin{align} \label{equation:cn<lcn<pcn}
\cn(G)\leq \lcn(G)\leq \pcn(G) \leq \De(G)+1,
\end{align}
where the last inequality follows easily from a greedy procedure.

The chromatic number $\cn$ is one of the most popular graph parameters. In 1941, Brooks \cite{Brooks41} proved that a connected graph $G$ satisfies $\cn(G)=\De(G)+1$ if and only if $G$ is a complete graph or an odd cycle. The list chromatic number $\lcn$ was introduced by 
Vizing \cite{Vizing76}. A few years later, Erd\H{o}s, Rubin, and Taylor \cite{ERT79} and Borodin \cite{BorodinThesis} independently proved that a connected graph $G$ satisfies $\lcn(G)=\De(G)+1$ if and only if $G$ is a complete graph or an odd cycle. The DP-chromatic number was introduced in 2018 by Dvo\v{r}\'ak and Postle \cite{DvorakP2018}; they used the term corresponding coloring. In 2018, Bernshteyn, Kostochka, and Pron \cite{BernshteynKP2017} proved that a connected graph $G$ satisfies $\pcn(G)=\De(G)+1$ if and only if $G$ is a 
complete graph or a cycle. Generalized DP-coloring of graphs were introduced by Kostochka, Schweser, and Stiebitz \cite{KostochkaSS2023}.

\medskip

Until now, generalized coloring of digraphs was only considered for the digraph properties $\cAD$ and $\cSD_m$ with $m\geq 1$. Note that the digraph property $\cSD_m$ with $m\geq 1$ is strongly reliable and $\cAD=\cSD_1$. Furthermore, it is easy to check that we have
\begin{align} \label{equation:CR(AD)}
\CR(\cAD)=\set{\dic{n}}{n\geq 2} \mbox{~~~and~~~} d(\cAD)=(1,1).
\end{align}
The only papers known to us that deal with $\cSD_m$-coloring of digraphs for arbitrary $m\geq 1$ are the papers by Golowich \cite{Golowich2016}, by Bang-Jensen \etal \cite{BangJensenSS2020}, and by Gon\c{c}alves \etal \cite{GoncalvesPR2024}. For a digraph $D$, let
$$\acn(D)=\acn(D:\cAD), \mbox{~~}\lacn(D)=\lacn(D:\cAD), \mbox{~~and~~} \pacn(D)=\pacn(D:\cAD);$$
the corresponding terms are respectively the \DF{dichromatic number}, the \DF{list-dichromatic number}, and the \DF{DP-dichromatic number} of $D$. The dichromatic number was introduced in 1982 by Neumann-Lara \cite{NeumannLara82}, but in the following two decades only a few papers appeared on the subject (see e.g. \cite{Erdos79, ErdosNL82, NeumannLara85, NeumannLara94}). The dichromatic number was rediscovered
at the beginning of the millennium by Mohar \cite{Mohar2003} and by Bokal, Fijav\v{z}, Juvan,
Kayll, and Mohar \cite{BokalFJKM2004}, and since then the topic has flourished and attracted
much more interest. The reader will find a long list of relevant references in \cite[Sections 3.10 and 5.8]{StiebitzST2024}. The list-dichromatic number was introduced by Harutyunyan and Mohar \cite{HarutMo2011}, and the DP-dichromatic number was introduced by Bang-Jensen, Bellitto, Schweser, and Stiebitz \cite{BangBSS2018}. Using a sequential coloring argument it is easy to show that every digraph $D$ satisfies
\begin{align} \label{equation:acn<lacn<pacn}
\acn(D)\leq \lacn(D)\leq \pacn(D)\leq \max \{\De^+(D),\De^-(D)\} + 1.
\end{align}

Bang-Jensen \etal \cite{BangBSS2018} proved that, if $D$ is connected, then $\pacn(D)=\max \{\De^+(D),\De^-(D)\} + 1$ if and only if $D$ is a directed cycle, a bidirected cycle, or a bidirected complete graph. This generalizes results obtained by Mohar \cite{Mohar2010} for the dichromatic number and by Harutyunyan and Mohar \cite{HarutMo2011} for the list-dichromatic number.

Let $G$ be a graph, let $D=D^\pm(G)$ be the corresponding bidirected graph, and let $\f:V(G)\to C$ be a coloring of $G$ and $D$, respectively. Evidently, $\f$ is an $\cO$-coloring of $G$ if and only if $\f$ is an $\cAD$-coloring of $D$. Consequently, every graph $G$ satisfies
\begin{align} \label{equation:D=G}
\acn(D^\pm(G))=\cn(G) \mbox{ and } \lacn(D^\pm(G))=\lcn(G).
\end{align}
Furthermore, as proved by Bang-Jensen \etal \cite{BangBSS2018}, every graph $G$ satisfies 
\begin{align} \label{equation:PD=PG}
\pacn(D^\pm(G))=\pcn(G). 
\end{align}
Note that \eqref{equation:PD=PG} is not obvious, since if $(X,H)$ is a cover of a bidirected graph $D$, then the digraph $H$ is not necessarily bidirected.  By \eqref{equation:D=G} and \eqref{equation:PD=PG}, results about $\cAD$-coloring of digraphs lead to results about $\cO$-coloring of graphs. 

Another interesting class of digraphs are oriented graphs. A digraph $D$ is called an \DF{orientation} of a graph $G$ if $D$ has no digons and $\UG(D)=G$. Erd\H{o}s and Neumann-Lara~\cite{Erdos79} conjectured that if a digraph $D$ is an orientation of a graph with maximum degree $\De$, then $\acn(D)= O(\De/\log \De)$. Neumann-Lara \cite{NeumannLara85} proposed the conjecture that if $D$ is an orientation of a planar graph, then $\acn(D)\leq 2$.
Both conjectures are widely open, see respectively \cite{Golowich2016,HarutyunyanEJC2011, PicasarriJGT2024, KawarabayashiSODA25, KawarabayashiArxiv25} and \cite{HarutMo2017,LiMohar2016}  for partial results.

\section{Criticality and main results}
\label{sect:critical+main}


In what follows, let $\cD$ be a reliable digraph property and $\pa$ be one of the $\cD$-coloring parameters $\acn, \lacn, \pacn$.
A digraph $D$ is called \DF{$(\cD,\pa)$-critical} and \DF{$(\cD,\pa,k)$-critical} if every proper induced subdigraph $D'$ of $D$ satisfies $\pa(D':\cD)<\pa(D:\cD)=k$, where $k\in \nato$. Furthermore, let
$$\Critd{k}=\set{D\in \cDG}{D \mbox{ is } (\cD,\pa,k)\mbox{-critical}}.$$
Proposition~\ref{prop:basicA}\ref{enum:prop:chi:0and1:a} immediately implies that 
\[
\Critd{0}=\{\ems\} \text{~~~and~~~} \Critd{1}=\{K_1\}.
\]
From Proposition~\ref{prop:basicA}\ref{enum:prop:chi:0and1:b}\ref{enum:prop:chi:0and1:c} it follows that
$$\Critd{2}=\CR(\cD)$$
unless $\pa=\pacn$ and $\cD$ is not strongly reliable. From Proposition~\ref{prop:basicB}\ref{prop:basic_properties_reliable:b}, it follows that every $(\cD,\pa)$-critical digraph is empty or connected.

\begin{proposition}
\label{prop:critsub}
Let $\cD$ be a reliable digraph property, let $D$ be a digraph, and let $\pa$ be one of the $\cD$-coloring parameters $\acn, \lacn, \pacn$. Then $D$ has an induced subgraph $D'$ such that $\pa(D':\cD)=\pa(D:\cD)$ and $D'$ is $(\cD,\pa)$-critical.
\end{proposition}
\begin{proof}
Among all subdigraphs $D'$ of $D$ satisfying $\pa(D')=\pa(D)$, we choose one with minimum order. Then $D'$ exists and has the desired properties by Proposition~\ref{prop:basicB}\ref{prop:basic_properties_reliable:a}. 
\end{proof}

The proposition above implies that many problems related to the $\cD$-coloring parameter $\pa$ can be reduced to problems about $(\cD,\pa)$-critical digraphs. The concept of critical graphs with respect to the chromatic number was introduced and first studied by G. A.~Dirac in the 1950s (see e.g. \cite{Dirac52, Dirac53, Dirac57, Dirac74}). Until now, the study of $\cn$-critical graphs 
has attracted a lot of attention and the reader will find a long list of relevant results in the book \cite{StiebitzST2024} ({\it Brooks' Theorem: Graph Coloring and Critical Graphs}). 

\medskip

Let $D$ be a digraph and $(X,H)$ be a cover of $D$. Given a vertex $v\in V(D)$, a partial transversal $T$ of $(X,H)$ such that $\dom(T:D)=V(D-v)$ and $H[T]\in \cD$ is said to be a \DF{$(\cD,v)$-transversal} of $(X,H)$. We call $(X,H)$ a \DF{$\cD$-critical cover} of $D$ if $(X,H)$ has no $\cD$-transversal, but for every vertex $v\in V(D)$ there exists a $(\cD,v)$-transversal.

\begin{proposition}
\label{prop:critcover}
Let $\cD$ be a reliable digraph property, let $D$ be a digraph, and let $k\geq 2$. Then the following statements hold:
\begin{enumerate}[label={\rm(\alph*)}]
\item If $D\in \Critc{k}$, then $D$ has a $\cD$-critical $(k-1)$-uniform cover.
\label{prop:critcover:a}
\item If $D\in \Critb{k}$, then $D$ has a $(k-1)$-assignment $L$ such that $C(D,L)$ is a $\cD$-critical $(k-1)$-uniform cover of $D$.
\label{prop:critcover:b}
\item If $D\in \Crita{k}$, then $C(D,L \equiv [1,k-1])$ is a $\cD$-critical $(k-1)$-uniform cover of $D$.
\label{prop:critcover:c}
\end{enumerate}
\end{proposition}
\begin{proof}
To prove~\ref{prop:critcover:a}, assume that $D\in \Critc{k}$. Since $\pdcn{D}> k-1$, there is a $(k-1)$-uniform cover $(X,H)$ of $D$ such that $(X,H)$ has no $\cD$-transversal. If $v\in V(D)$, then $\pdcn{D-v}\leq k-1$ implies that $(X,H)$ has a $(\cD,v)$-transversal. Hence $(X,H)$ is a $\cD$-critical $(k-1)$-uniform cover of $D$. This proves~\ref{prop:critcover:a}. The proofs for~\ref{prop:critcover:b} and~\ref{prop:critcover:c} are similar.
\end{proof}

\begin{proposition}
\label{prop:lowvertex}
Let $\cD$ be a reliable digraph property, let $D$ be a digraph, and let $(X,H)$ be a $\cD$-critical cover of $D$.
Then the following statements hold.
\begin{enumerate}[label={\rm (\alph*)}]
\item  For every $v\in V(D)$, $d_D(v) \geq d(\cD) \cdot |X_v|$.
\label{enum:lowvertex:a}
\item Let $v\in V(D)$ be a vertex such that $d^+_D(v)=d^+(\cD)|X_v|$, and let $T$ be a $(\cD,v)$-transversal of $(X,H)$. Moreover, for $x\in X_v$, let $H_x=H[T \cup \{x\}]$ and $d_x^+=d^+_{H_x}(x)$. Then $d^+_x=d^+(\cD)$ for all $x\in X_v$ and $d_D^+(v)=\sum_{x\in X_v}d^+_x$.
\label{enum:lowvertex:b}
\item Let $v\in V(D)$ be a vertex such that $d^-_D(v)=d^-(\cD)|X_v|$, and let $T$ be a $(\cD,v)$-transversal of $(X,H)$. Moreover, for $x\in X_v$, let $H_x=H[T \cup \{x\}]$ and $d_x^-=d^-_{H_x}(x)$. Then $d^-_x=d^-(\cD)$ for all $x\in X_v$ and $d_D^-(v)=\sum_{x\in X_v}d^-_x$.
\label{enum:lowvertex:c}
\end{enumerate}
\end{proposition}
\begin{proof}
We prove the proposition for the out-degree; the proof for the in-degree is similar. Let $v\in V(D)$ be an arbitrary vertex. Since $(X,H)$ is a $\cD$-critical cover of $D$, there is a $(\cD,v)$-transversal $T$ of $(X,H)$. Since $(X,H)$ has no $\cD$-transversal, for every color $x\in X_v$, $H_x=H[T \cup \{x\}]\not\in \cD$. Then Proposition~\ref{prop:smooth}\ref{enum:prop:smooth:e} implies that $d_x^+=d^+_{H_x}(x)\geq d^+(\cD)$ for all $x\in X_v$. For every vertex $u\in V(D-v)$, there is a unique vertex $x_u$ with $T\cap X_u=\{x_u\}$. From Proposition~\ref{prop:cover} we then obtain that
\begin{eqnarray*}
d_D^+(v) &= & \sum_{u\in V(D-v)} a_D(v,u) \geq \sum_{u\in V(D-v)} \sum_{x\in X_v}a_H(x,x_u)\\
       &\geq & \sum_{x\in X_v}\sum_{u\in V(D-v)}a_H(x,x_u)=\sum_{x\in X_v}d^+_x\geq d^+(\cD)|X_v|.
\end{eqnarray*}
Then $d_D^+(v)=d^+(\cD)|X_v|$ implies that $d_x^+=d^+(\cD)$ for all $x\in X_v$. This proves~\ref{enum:lowvertex:a} for the out-degree and~\ref{enum:lowvertex:b}.
\end{proof}

Let $\cD$ be a reliable digraph property, let $D$ be a digraph, and let $(X,H)$ be a $\cD$-critical cover of $D$. Then, we define
$$V(D,X,H,\cD)=\set{v\in V(D)}{d_D(v)=d(\cD)\cdot |X_v|}.$$
A vertex $v\in V(D)$ is called a \DF{low vertex} of $D$ with respect to $(X,H,\cD)$ if $v\in V(D,X,H,\cD)$, and a \DF{high vertex} otherwise. Moreover, we call $D[V(D,X,H,\cD)]$ the \DF{low vertex subdigraph} of $D$ with respect to $(X,H,\cD)$.

A digraph $D$ is called a \DF{$\cD$-dibrick} if $D$ is a connected digraph such that $D=D^\pm(K_n)$ with $n\geq 1$, or $D=D^\pm(C_n)$ with odd $n\geq 3$, or $D \in \cD$  with
$\De(D)\leq d(\cD)$, or $D\in \CR(\cD)$ is an Eulerian digraph such that $d_D^+(v)=d^-_D(v)=d^+(\cD)=d^-(\cD)$ for all $v\in V(D)$.

\smallskip

The idea of dividing the vertices of a $\cn$-critical graph into low and high vertices dates back to the seminal 1963 papers \cite{Gallai63a} and \cite{Gallai63b} by T. Gallai. If $G$ is a $\cn$-critical graph with $\cn(G)=k\geq 1$ (i.e. $\cn(G')<k$ for every induced subgraph $G'$ of $G$), then $\de(G)\geq k-1$ and the low vertices of $G$ are the vertices having degree $k-1$ in $G$. Gallai \cite{Gallai63a} proved that the low vertex subgraph of a $\cn$-critical graph has a simple block structure, each of its blocks being a complete graph or an odd cycle. That the same statement holds for $\lcn$-critical graphs was proved by Erd\H{o}s, Rubin and Taylor \cite{ERT79}. For the class of $\pcn$-critical graphs, their low vertex subgraph can have an additional type of blocks, namely even cycles, as proved by Bernshteyn, Kostochka and Pron \cite{BernshteynKP2017}.  Using \eqref{equation:D=G} and \eqref{equation:PD=PG}, it is easy to show that all three results about the block structure of critical graphs can be obtained from the following two theorems, by considering the digraph property $\cD=\cAD$. Thus, the following two theorems can be seen as a far-reaching extension of Gallai's fundamental result. 
The case of $L$ being constant in the second theorem has been proved by Bang-Jensen, Bellitto, Stiebitz, and Schweser~\cite{BangJensenBSS2019}.

\begin{restatable}{theorem}{theoremmainDP}
\label{theorem:mainDP}
Let $\cD$ be a reliable digraph property, let $D$ be a digraph, and let $(X,H)$ be a $\cD$-critical cover of $D$. If $B$ is a block of the low vertex subdigraph $D[V(D,X,H,\cD)]$, then $B$ is a $\cD$-dibrick, or a bidirected cycle of even order, or an antidirected cycle.
\end{restatable}

\begin{restatable}{theorem}{theoremmainLI}
\label{theorem:mainLI}
Let $\cD$ be a reliable digraph property, let $D$ be a digraph, let $L$ be a list assignment, and let $(X,H)=C(D,L)$ be a $\cD$-critical cover of $D$. If $B$ is a block of the low vertex subdigraph $D[V(D,X,H,\cD)]$, then $B$ is a $\cD$-dibrick. 
\end{restatable}

The proofs of the two theorems above are postponed to Section~\ref{sectsub:proofs:mainDP_LI} because they rely on our third main result presented in Section~\ref{sect:DP-col+degeneracy}. The rest of this section discusses some consequences of the two theorems.  

\begin{corollary}
\label{corollary:degree-choosable}
Let $\cD$ be a reliable digraph property, let $D$ be a connected digraph, $|D| \geq 2$,
and let $(X,H)$ be a cover of $D$ such that
each vertex $v\in V(D)$ satisfies $d(v) \leq d(\cD)\cdot |X_v|$. If $D$ has no $(\cD,(X,H))$-coloring, then $d^+(\cD)=d^-(\cD)=r\geq 1$, and the following statements hold.
\begin{enumerate}[label={\rm (\alph*)}]
\item Every block $B$ of $D$ satisfies $B=D^\pm(K_n)$ with $n\geq 2$, or $B=D^\pm(C_n)$ with $n\geq 3$, or $B\in \cD\cup \CR(\cD)$ is $r$-diregular.
\label{enum:cor:degree-choosable:a}
\item If $(X,H)$ is associated with a list assignment of $D$, then every block $B$ of $D$ satisfies $B=D^\pm(K_n)$ with $n\geq 2$, or $B=D^\pm(C_n)$ with odd $n\geq 3$, or $B\in \cD \cup \CR(\cD)$ is $r$-diregular.
\label{enum:cor:degree-choosable:b}
\end{enumerate}
\end{corollary} 
\begin{proof}
Since $D$ has no $(\cD,(X,H))$-coloring, there exists an induced subdigraph $D'$ of $D$ such that $(X,H)/D'$ is a $\cD$-critical cover of $D'$. Since $\cD$ is additive, it follows that $D'$ is connected.
Using Proposition~\ref{prop:lowvertex}\ref{enum:lowvertex:a}, every vertex $u\in V(D')$ satisfies
$$d_{D'}(u)\geq d(\cD) \cdot |X_u|\geq d_D(u) .$$
Since $D$ is connected and $D'\not=\ems$, this implies that $D'=D$ and, every vertex $u\in V(D)$ satisfies
$$d_{D}(u)= d(\cD)\cdot |X_u|.$$
Consequently, $(X,H)$ is a $\cD$-critical cover of $D$ and each vertex of $D$ is a low vertex. Since $|D|\geq 2$, we obtain that
$$0<|A(D)|=\sum_{u\in V(D)} d^+_D(u)=d^+(\cD) \sum_{u\in V(D)}|X_u|=\sum_{u\in V(D)} d^-_D(u)=d^-(\cD) \sum_{u\in V(D)}|X_u|.$$
This implies that $d^+(\cD)=d^-(\cD)=r\geq 1$ and $d_D^+(u)=d^-_D(u)=r|X_u|$ for all $u\in V(D)$. In particular, $D$ is an Eulerian digraph. Using induction on the number of blocks of $D$, it then follows that each block of $D$ is an Eulerian digraph too.  
Hence~\ref{enum:cor:degree-choosable:a} follows from Theorem~\ref{theorem:mainDP} and~\ref{enum:cor:degree-choosable:b} follows from Theorem~\ref{theorem:mainLI}.
\end{proof}

\begin{corollary}
\label{corollary:brooksDP}
Let $\cD$ be a reliable digraph property such that $d(\cD)\geq (1,1)$, and let $k\in \nato$. Then the following statements hold:
\begin{enumerate}[label={\rm (\alph*)}]
\item If $D\in \Critc{k+1}$, then $\de(D) \geq d(\cD) \cdot k$. Moreover, if $U=\set{v\in V(D)}{d_D(v)=d(\cD)\cdot k}$, then every block of $D[U]$ is a $\cD$-dibrick, or a bidirected cycle of even order, or an antidirected cycle.
\label{enum:cor:brooksDP:a}
\item If $G\in \Crita{k+1}\cup \Critb{k+1}$, then $\de(D)\geq d(\cD)\cdot k$. Moreover, if $U=\set{v\in V(D)}{d_D(v)=d(\cD)\cdot k}$, then every block of $D[U]$ is a $\cD$-dibrick. 
\label{enum:cor:brooksDP:b}
\item Every digraph $D$ satisfies $\displaystyle\pdcn{D}\leq \min \bigg\{\frac{\De^+(D)}{d^+(\cD)}, \frac{\De^-(D)}{d^-(\cD)} \bigg\}+1$.
\label{enum:cor:brooksDP:c}
\end{enumerate}
\end{corollary}
\begin{proof}
Statements~\ref{enum:cor:brooksDP:a} and~\ref{enum:cor:brooksDP:b} are evident if $k=0$, since 
$$\Crita{1}=\Critb{1}=\Critc{1}=\{K_1\}.$$
To prove~\ref{enum:cor:brooksDP:a}, assume that $D\in \Critc{k+1}$ with $k\geq 1$. Then there is a $\cD$-critical $k$-uniform cover of $D$ (by Proposition~\ref{prop:critcover}). Hence, every vertex $v\in V(D)$ satisfies that $d_D(v)\geq d(\cD)\cdot k$ (by Proposition~\ref{prop:lowvertex}), and thus $\de(D)\geq d(\cD)\cdot k$. Furthermore, we have $V(D,X,H,\cD)=U$. Therefore, the statements about the block structure of $D[U]$ follow from Theorem~\ref{theorem:mainDP}. This proves~\ref{enum:cor:brooksDP:a}. The proof of~\ref{enum:cor:brooksDP:b} is similar using Theorem~\ref{theorem:mainLI}.

For the proof of~\ref{enum:cor:brooksDP:c}, let $D$ be an arbitrary digraph with $\pdcn{D}=k+1$ and $k\geq 0$. Then there exists an induced subdigraph $D'$ of $D$ with $D'\in \Critc{k+1}$ (by Proposition~\ref{prop:critsub}). Using~\ref{enum:cor:brooksDP:a}, we then obtain that
$$k\cdot d(\cD) \leq \de(D')\leq \De(D),$$
which gives $k\leq \min \{\De^+(D)/d^+(\cD), \De^-(\cD)/d^-(\cD)\}$. This proves~\ref{enum:cor:brooksDP:c}.
\end{proof}

We conclude the section with two  Brooks' type results for the $\cD$-coloring parameters. The proof of the second theorem is similar to the proof of the first theorem and we omit it.

\begin{theorem}
\label{theorem:brooksPD}
Let $\cD$ be a strongly reliable digraph property such that $d(\cD)\geq (1,1)$, and let $D$ be a connected digraph. Then 
$$\pdcn{D}\leq \max \bigg\{ \left\lceil \frac{\De^+(D)}{d^+(\cD)} \right\rceil, \left\lceil\frac{\De^-(D)}{d^-(\cD)} \right\rceil \bigg\}$$
unless $D=K_1$, or $d^+(\cD)=d^-(\cD)=r\geq 1$ and one of the following holds:
\begin{itemize}
	\item $D\in \CR(\cD)=\Critc{2}$ and $D$ is $r$-diregular, 
	\item $D=D^\pm(K_{kr+1})\in \Critc{k+1}$ with $k\geq 1$, or
	\item $r=1$ and $D=D^\pm(C_n)\in \Critc{3}$ with $n\geq 3$.
\end{itemize}
\end{theorem}
\begin{proof}
Let $D$ be a connected digraph. To prove the upper bound for $\pdcn{D}$, it is sufficient to consider the case where $k=\De^+(D)/d^+(\cD)=\De^-(D)/d^-(\cD)$ is an integer and $\pdcn{D}=k+1$, otherwise the result follows from Corollary~\ref{corollary:brooksDP}\ref{enum:cor:brooksDP:c}. Then $D$ contains an induced subdigraph $D'$ such that $D'\in \Critc{k+1}$ (by Proposition~\ref{prop:critsub}).

By Corollary~\ref{corollary:brooksDP}\ref{enum:cor:brooksDP:a}, we then have $\de(D')\geq k\cdot  d(\cD)=\De(D)$ . Since $D'$ is an induced subdigraph of $D$ and $D$ is connected, this implies that $D=D'\in \Critc{k+1}$ and, for every vertex $v\in V(D)$, we have $d_D(v)=k \cdot d(\cD)$. If $k=0$, then $D=K_1$ and we are done, so assume that $k\geq 1$. 

We have $0<|A(D)|=kd^+(\cD)|D|=kd^-(\cD)|D|$, which gives $d^+(\cD)=d^-(\cD)=r$ with $r\in \nat$. Furthermore, for every vertex $v\in V(D)$, we have $d^+_D(v)=d^-_D(v)=kr$. Hence every vertex of $D$ is a low vertex, and every block of $D$ is a $\cD$-dibrick, or a bidirected even cycle, or an antidirected cycle (by Corollary~\ref{corollary:brooksDP}\ref{enum:cor:brooksDP:a}). 

Assume first that $k=1$. Then $D$ is $r$-diregular and, since $\cD$ is strongly reliable, $D\in \Critc{2}=\CR(\cD)$ by Proposition, so we are done. Henceforth, we assume that $k\geq 2$.

It is well-known and easy to prove by induction that every block of an Eulerian digraph is Eulerian itself. Therefore, every block $B$ of $D$ satisfies $B=D^\pm(K_n)$ with $n\geq 2$, $B=D^\pm(C_n)$ with $n\geq 2$, or $B\in \cD\cup \CR(\cD)$ is $r$-diregular. Since $k\geq 2$ and $D$ is $kr$-diregular, it is straightforward to prove that $D$ is a block, and either $D=D^\pm(K_n)$ with $n\geq 2$ or $D=D^\pm(C_n)$ with $n\geq 3$.

If $D=D^\pm(K_n)$, then we get $kr=d_D^+(v)=d_D^-(v)=n-1$ and so $n=kr+1$; hence we are done. If $D=D^\pm(C_n)$, then we get $kr=d_D^+(v)=d_D^-(v)=2$. Since $k\geq 2$, this leads to $k=2$, $r=1$, and so $D\in \Critc{3}$. This completes the proof.  
\end{proof}

\begin{theorem}
\label{theorem:brooksLI}
Let $\cD$ be a reliable digraph property such that $d(\cD)\geq (1,1)$, and let $D$ be a connected digraph. Then 
$$\ldcn{D}\leq \max \bigg\{ \left\lceil \frac{\De^+(D)}{d^+(\cD)} \right\rceil, \left\lceil\frac{\De^-(D)}{d^-(\cD)} \right\rceil \bigg\}$$
unless $D=K_1$, or $d^+(\cD)=d^-(\cD)=r\geq 1$ and one of the following holds:
\begin{itemize}
	\item $D\in \CR(\cD)=\Critb{2}$ and $D$ is $r$-diregular, 
	\item $D=D^\pm(K_{kr+1})\in \Critb{k+1}$ with $k\geq 1$, or 
	\item $r=1$ and $D=D^\pm(C_n)\in \Critb{3}$ with odd $n\geq 3$.
\end{itemize} 
\end{theorem}

\section{DP-coloring and variable degeneracy}
\label{sect:DP-col+degeneracy}

In this section we shall establish a result (Theorem~\ref{main_theorem}) that combines DP-coloring with variable degeneracy. 

\subsection{Variable degeneracy and constructible configurations}

Let $H$ be a digraph and let $f$ be a \DF{vertex function} of $H$, that is, $f:V(H)\to \nato^2$. We denote by $f^+$ and $f^-$ the projections of $f$ on the first and second coordinates respectively; so $f(v)=(f^+(v),f^-(v))$ for $v\in V(D)$.
Furthermore, $\su(f)=\set{x\in V(H)}{f(x)\not=(0,0)}$ is the \DF{support} of $f$ in $H$, and $\suc(f)=\set{x\in V(H)}{f(x)=(0,0)}$ is the \DF{complementary support}
of $f$ in $H$. For a function $g:V(H)\to \nato$ (or $g:V(H)\to \nato^2$) and a set $X\subseteq V(H)$, define
$$g(X)=\sum_{x\in X}g(x).$$
We say that $H$ is \DF{strictly $f$-degenerate} if each nonempty subdigraph $H'$ of $H$ contains a vertex (color) $x$ such that $d_{H'}^+(x)<f^+(x)$ or $d_{H'}^-(x)<f^-(x)$. Note that if a subdigraph $\tilde{H}$ of $H$ is a strictly $f$-degenerate, then $V(\tilde{H})\subseteq \su(f)$. The concept of variable degeneracy for graphs seems to have been first studied by Borodin, Kostochka, and Toft \cite{BorodinKT00}. This concept was extended to digraphs by Bang-Jensen, Schweser, and Stiebitz \cite{BangJensenSS2020} for the symmetric case, where $f^+=f^-$, and by Gon\c{c}alves, Picasarri-Arrieta, and Reinald \cite{GoncalvesPR2024} for the general case. 

\smallskip

In this section we deal with the following coloring problem. A \DF{configuration} is a tuple $\cK=
(D,X,H,f)$ such that $D$ is a digraph, $(X,H)$ is a cover of $D$, and $f$ is a vertex function of $H$. Given a configuration $\cK=(D,X,H,f)$, we want to decide whether $(X,H)$ has a transversal $T$ such that $H[T]$ is strictly $f$-degenerate. In general, this decision problem is {\NPC}. However, if we add a certain degree condition it might become a polynomial problem.

In what follows, let $\cK=(D,X,H,f)$ be a configuration. We call $\cK$  \DF{degree-feasible} if for each vertex $v$ of $D$ we have
$$f(X_v)=\sum_{x \in X_v} f(x) \geq d_D(v).$$
Furthermore, we say that $\cK$ is \DF{colorable} if $(X,H)$ has a transversal $T$ such that $H[T]$ is strictly $f$-degenerate, otherwise $\cK$ is said to be \DF{uncolorable}.

In what follows, we give an exact characterization of degree-feasible uncolorable configurations.
For this, let us recursively define the family of \DF{constructible configurations}. A configuration $(D,X,H,f)$ is called \DF{constructible} if one of the following six conditions hold.

\begin{enumerate}[label=$\mathrm{(\Roman*)}$]
\item $(D,X,H,f)$ is an \upshape{M}-\DF{configuration}, that is, $D$ is a block and there exists a transversal $T$ of $(X,H)$ such that $(X,H)/T$ is a saturated cover of $D$, and for $v\in V(D)$ and $x\in X_v$, we have $f(x)=d_D(v)$ if $x\in T$ and $f(x)=(0,0)$ otherwise. See Figure~\ref{fig:m_configuration} for an illustration.

\begin{figure}
	\begin{center}
		\begin{tikzpicture}[thin,scale=1, every node/.style={transform shape}]
			\tikzset{vertex/.style = {circle,fill=black,minimum size=5pt, inner sep=0pt}}
			\tikzset{edge/.style = {->,> = latex}}
			
			\begin{scope}
				\node[vertex,label=below:$v_1$] (v1) at (-2, 0) {};
				\node[vertex,label=below:$v_2$] (v2) at (-0.66, 0) {};
				\node[vertex,label=below:$v_3$] (v3) at (0.66, 0) {};
				\node[vertex,label=below:$v_4$] (v4) at (2, 0) {};
				\node[] (D) at (-3.5, 0) {$D:$};
				\node[opacity=0] (D) at (3.5, 0) {$D:$};
				\draw[edge] (v1) to (v2);
				\draw[edge, bend left=30] (v1) to (v3);
				\draw[edge] (v3) to (v2);
				\draw[edge] (v4) to (v3);
				\draw[edge, bend right=30] (v4) to (v2);
			\end{scope}
			
			\begin{scope}[yshift=-1.2cm]
				\node[vertex,label=below:{\footnotesize $(2,0)$}] (v1) at (-2, 0) {};
				\node[vertex,label=below:{\footnotesize $(0,3)$}] (v2) at (-0.66, 0) {};
				\node[vertex,label=below:{\footnotesize $(1,2)$}] (v3) at (0.66, 0) {};
				\node[vertex,label=below:{\footnotesize $(2,0)$}] (v4) at (2, 0) {};
				\draw[edge] (v1) to (v2);
				\draw[edge, bend left=30] (v1) to (v3);
				\draw[edge] (v3) to (v2);
				\draw[edge] (v4) to (v3);
				\draw[edge, bend right=30] (v4) to (v2);
				\node[] (H) at (-3.5, 0) {$H:$};
				
				\draw[dotted] (-2,-0.15) ellipse (0.5cm and 0.7cm);
				\node[] (X1) at (-2,-1.2) {$X_{v_1}$};
				\draw[dotted] (-0.66,-0.15) ellipse (0.5cm and 0.7cm);
				\node[] (X1) at (-0.66,-1.2) {$X_{v_2}$};
				\draw[dotted] (0.66,-0.15) ellipse (0.5cm and 0.7cm);
				\node[] (X1) at (0.66,-1.2) {$X_{v_3}$};
				\draw[dotted] (2,-0.15) ellipse (0.5cm and 0.7cm);
				\node[] (X1) at (2,-1.2) {$X_{v_4}$};
			\end{scope}
			
		\end{tikzpicture}
		\caption{An illustration of an M-configuration $(D,X,H,f)$.}
		\label{fig:m_configuration}
	\end{center}
\end{figure}
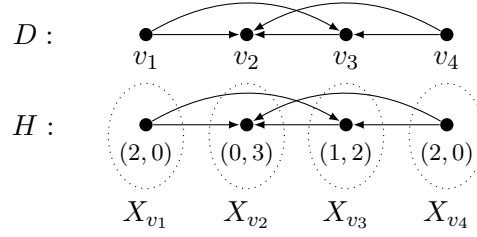

\item$(D,X,H,f)$ is a \upshape{K}-\DF{configuration}, that is, $D=D^\pm(K_n)$ with $n\in \nat$ and there are integers $n_1,n_2,\ldots,n_p\in \nat$ with $p \geq 1$ such that $n_1 + n_2 + \cdots + n_p=n-1$. Moreover, there are $p$ pairwise disjoint transversals $T_1, T_2, \ldots, T_p$ of $(X,H)$ such that $(X,H)/T_i$ is a saturated cover of $D$ for $i\in [1,p]$, and, for $v\in V(H)$ and $x\in X_v$, we have $f(x)=(n_i,n_i)$ if $x\in T_i$ for $i\in [1,p]$ and $f(x)=(0,0)$ otherwise. See Figure~\ref{fig:k_configuration} for an illustration.

\item $(D,X,H,f)$ is an \DF{odd} \upshape{C}-\DF{configuration}, that is, $D=D^\pm(C_n)$ with $n\in \nat$ and $n\geq 5$ odd. Moreover, there are two disjoint transversals $T_1, T_2$ of $(X,H)$ such that $(X,H)/T_i$ is a saturated cover of $D$ for $i\in\{1,2\}$, and,
    for $v\in V(D)$ and $x\in X_v$, we have $f(x)=(1,1)$ if $x\in T_1 \cup T_2$ and $f(x)=(0,0)$ otherwise. See Figure~\ref{fig:oddc_configuration} for an illustration.
    
    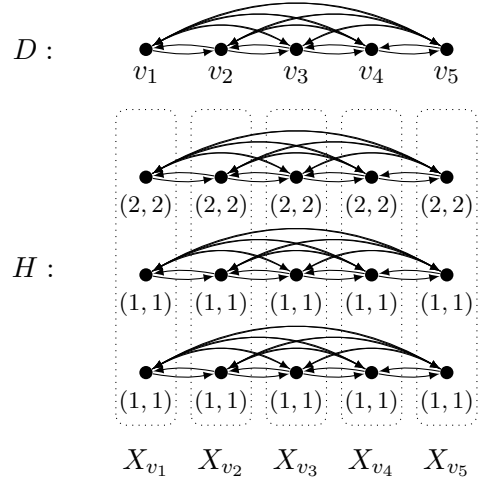
\begin{figure}
	\begin{center}
		\begin{tikzpicture}[thin,scale=1, every node/.style={transform shape}]
			\tikzset{vertex/.style = {circle,fill=black,minimum size=5pt, inner sep=0pt}}
			\tikzset{edge/.style = {->,> = latex}}
			
			\begin{scope}
				\node[] (D) at (-3.5, 0) {$D:$};
				\node[opacity=0] (D) at (3.5, 0) {$D:$};
				\foreach \i in {1,...,5}{
					\node[vertex,label=below:$v_{\i}$] (v\i) at (-3+\i, 0) {};
					\draw[rounded corners, dotted] (-3.4+\i, -5) rectangle (-2.6+\i, -0.8) {};
					\node[] (X1) at (-3+\i,-5.5) {$X_{v_\i}$};
					\ifthenelse{\i>1}{
						\pgfmathtruncatemacro{\j}{\i-1}
						\draw[edge, bend right=10] (v\i) to (v\j) {};
						\draw[edge, bend right=10] (v\j) to (v\i) {};
						\ifthenelse{\i>2}{
							\pgfmathtruncatemacro{\l}{\j-1}
							\foreach \k in {1,...,\l}{
								\draw[edge,bend right=30] (v\i) to (v\k) {};
								\draw[edge,bend left=30] (v\k) to (v\i) {};
							}
						}{}
					}{}
				}
			\end{scope}
			
			\foreach \a in {0,1,2}{
				\pgfmathsetmacro{\shift}{-1.7-1.3*\a}
			\begin{scope}[yshift=\shift cm]
				\foreach \i in {1,...,5}{
					\ifthenelse{\a=0}{
						\node[vertex,label=below:{\footnotesize $(2,2)$}] (v\i) at (-3+\i, 0) {};
					}
					{
						\node[vertex,label=below:{\footnotesize $(1,1)$}] (v\i) at (-3+\i, 0) {};
					}
					\ifthenelse{\i>1}{
						\pgfmathtruncatemacro{\j}{\i-1}
						\draw[edge, bend right=10] (v\i) to (v\j) {};
						\draw[edge, bend right=10] (v\j) to (v\i) {};
						\ifthenelse{\i>2}{
							\pgfmathtruncatemacro{\l}{\j-1}
							\foreach \k in {1,...,\l}{
								\draw[edge,bend right=30] (v\i) to (v\k) {};
								\draw[edge,bend left=30] (v\k) to (v\i) {};
							}
						}{}
					}{}
				}
			\end{scope}
			}
			\node[] (H) at (-3.5, -2.9) {$H:$};
		\end{tikzpicture}
		\caption{An illustration of a K-configuration $(D,X,H,f)$.}
		\label{fig:k_configuration}
	\end{center}
\end{figure}

    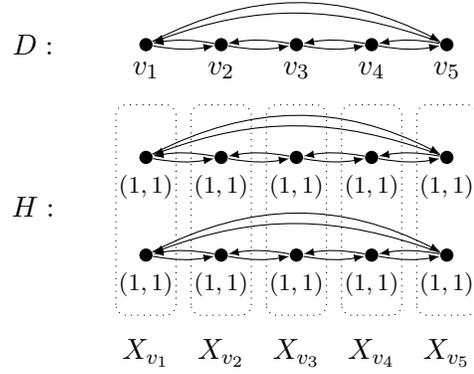
\begin{figure}
    	\begin{center}
    		\begin{tikzpicture}[thin,scale=1, every node/.style={transform shape}]
    			\tikzset{vertex/.style = {circle,fill=black,minimum size=5pt, inner sep=0pt}}
    			\tikzset{edge/.style = {->,> = latex}}
    			
    			\begin{scope}
    				\node[] (D) at (-3.5, 0) {$D:$};
    				\node[opacity=0] (D) at (3.5, 0) {$D:$};
    				\foreach \i in {1,...,5}{
    					\node[vertex,label=below:$v_{\i}$] (v\i) at (-3+\i, 0) {};
    					\draw[rounded corners, dotted] (-3.4+\i, -3.6) rectangle (-2.6+\i, -0.8) {};
    					\node[] (X1) at (-3+\i,-4.1) {$X_{v_\i}$};
    					\ifthenelse{\i>1}{
    						\pgfmathtruncatemacro{\j}{\i-1}
    						\draw[edge, bend right=10] (v\i) to (v\j) {};
    						\draw[edge, bend right=10] (v\j) to (v\i) {};
    					}{}
    				}
    				\draw[edge, bend right=20] (v5) to (v1) {};
    				\draw[edge, bend left=27] (v1) to (v5) {};
    			\end{scope}
    			
    			\foreach \a in {0,1}{
    				\pgfmathsetmacro{\shift}{-1.5-1.3*\a}
    				\begin{scope}[yshift=\shift cm]
    					\foreach \i in {1,...,5}{
    						\node[vertex,label=below:{\footnotesize $(1,1)$}] (v\i) at (-3+\i, 0) {};
    						\ifthenelse{\i>1}{
    							\pgfmathtruncatemacro{\j}{\i-1}
    							\draw[edge, bend right=10] (v\i) to (v\j) {};
    							\draw[edge, bend right=10] (v\j) to (v\i) {};
    						}{}
    					}
  						\draw[edge, bend right=20] (v5) to (v1) {};
  						\draw[edge, bend left=27] (v1) to (v5) {};
    				\end{scope}
    			}
    			\node[] (H) at (-3.5, -2.15) {$H:$};
    		\end{tikzpicture}
    		\caption{An illustration of an odd C-configuration $(D,X,H,f)$.}
    		\label{fig:oddc_configuration}
    	\end{center}
    \end{figure}

\item $(D,X,H,f)$ is an \DF{even} \upshape{C}-\DF{configuration}, that is, $D=D^\pm(C_n)$ with $n\in \nat$ and $n\geq 4$ even. Moreover, there are two disjoint transversals $T_1, T_2$ of $(X,H)$ such that $H[T_1\cup T_2] = D^\pm(C_{2n})$, and for $v\in V(D)$ and $x\in X_v$, we have $f(x)=(1,1)$ if $x\in T_1\cup T_2$ and $f(x)=(0,0)$ otherwise. See Figure~\ref{fig:evenc_configuration} for an illustration.

\item $(D,X,H,f)$ is an \upshape{A}-\DF{configuration}, that is, $D$ is the antidirected cycle on $n$ vertices with $n\in \nat$ and $n\geq 4$ even. Moreover, there are two disjoint transversals $T_1, T_2$ of $(X,H)$ such that $H[T_1\cup T_2]$ is the antidirected cycle on $2n$ vertices, and for $v\in V(D)$ and $x\in X_v$, we have
$$f(x)=\begin{cases}
	(1,0) & \text{if } x\in T_1\cup T_2 \text{ and } v \text{ is a source of }$D$,\\
	(0,1) & \text{if } x\in T_1\cup T_2 \text{ and } v \text{ is a sink of }$D$,\\
	(0,0) & \text{otherwise.}
\end{cases}$$
See Figure~\ref{fig:a_configuration} for an illustration.

\begin{figure}
	\begin{center}
		\begin{tikzpicture}[thin,scale=1, every node/.style={transform shape}]
			\tikzset{vertex/.style = {circle,fill=black,minimum size=5pt, inner sep=0pt}}
			\tikzset{edge/.style = {->,> = latex}}
			
			\begin{scope}
				\node[] (D) at (-4, 0) {$D:$};
				\node[opacity=0] (D) at (4, 0) {$D:$};
				\foreach \i in {1,...,6}{
					\node[vertex,label=below:$v_{\i}$] (v\i) at (-3.5+\i, 0) {};
					\draw[rounded corners, dotted] (-3.9+\i, -3.6) rectangle (-3.1+\i, -0.8) {};
					\node[] (X1) at (-3.5+\i,-4.1) {$X_{v_\i}$};
					\ifthenelse{\i>1}{
						\pgfmathtruncatemacro{\j}{\i-1}
						\draw[edge, bend right=10] (v\i) to (v\j) {};
						\draw[edge, bend right=10] (v\j) to (v\i) {};
					}{}
				}
				\draw[edge, bend right=20] (v6) to (v1) {};
				\draw[edge, bend left=27] (v1) to (v6) {};
			\end{scope}
			
			\foreach \a in {0,1}{
				\pgfmathsetmacro{\shift}{-1.5-1.3*\a}
				\begin{scope}[yshift=\shift cm]
					\foreach \i in {1,...,6}{
						\ifthenelse{\a=0}{
							\node[vertex,label=above:{\footnotesize $(1,1)$}] (v\i\a) at (-3.5+\i, 0) {};
						}{
							\node[vertex,label=below:{\footnotesize $(1,1)$}] (v\i\a) at (-3.5+\i, 0) {};
						}
						\ifthenelse{\i>1}{
							\pgfmathtruncatemacro{\j}{\i-1}
							\draw[edge, bend right=10] (v\i\a) to (v\j\a) {};
							\draw[edge, bend right=10] (v\j\a) to (v\i\a) {};
						}{}
					}
				\end{scope}
			}
			\draw[edge, bend right=3] (v60) to (v11) {};
			\draw[edge, bend right=3] (v11) to (v60) {};
			\draw[edge, bend right=3] (v61) to (v10) {};
			\draw[edge, bend right=3] (v10) to (v61) {};
			\node[] (H) at (-4, -2.15) {$H:$};
		\end{tikzpicture}
		\caption{An illustration of an even C-configuration $(D,X,H,f)$.}
		\label{fig:evenc_configuration}
	\end{center}
\end{figure}
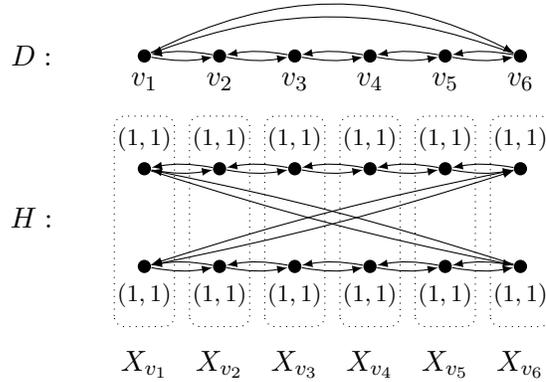

\begin{figure}
	\begin{center}
		\begin{tikzpicture}[thin,scale=1, every node/.style={transform shape}]
			\tikzset{vertex/.style = {circle,fill=black,minimum size=5pt, inner sep=0pt}}
			\tikzset{edge/.style = {->,> = latex}}
			
			\begin{scope}
				\node[] (D) at (-4, 0) {$D:$};
				\node[opacity=0] (D) at (4, 0) {$D:$};
				\foreach \i in {1,...,6}{
					\draw[rounded corners, dotted] (-3.9+\i, -3.6) rectangle (-3.1+\i, -0.8) {};
					\node[] (X1) at (-3.5+\i,-4.1) {$X_{v_\i}$};
				}
				\foreach \i in {1,3,5}{
					\node[vertex,label=below:{$v_\i$}] (v\i) at (-3.5+\i, 0) {};
				}
				\foreach \i in {2,4,6}{
					\node[vertex,label=below:{$v_\i$}] (v\i) at (-3.5+\i, 0) {};
					\pgfmathsetmacro{\j}{\i-1}
					\draw[edge] (v\j) to (v\i) {};
					\ifthenelse{\i<6}{
						\pgfmathsetmacro{\k}{\i+1}
						\draw[edge] (v\k) to (v\i) {};
					}{}
				}
				\draw[edge, bend left=20] (v1) to (v6);
			\end{scope}
			
			\foreach \a in {0,1}{
				\pgfmathsetmacro{\shift}{-1.5-1.3*\a}
				\begin{scope}[yshift=\shift cm]
					\foreach \i in {1,3,5}{
						\ifthenelse{\a=0}{
							\node[vertex,label=above:{\footnotesize $(1,0)$}] (u\i\a) at (-3.5+\i, 0) {};
						}{
							\node[vertex,label=below:{\footnotesize $(1,0)$}] (u\i\a) at (-3.5+\i, 0) {};
						}
					}
					\foreach \i in {2,4,6}{
						\ifthenelse{\a=0}{
							\node[vertex,label=above:{\footnotesize $(0,1)$}] (u\i\a) at (-3.5+\i, 0) {};
						}{
							\node[vertex,label=below:{\footnotesize $(0,1)$}] (u\i\a) at (-3.5+\i, 0) {};
						}
						\pgfmathtruncatemacro{\j}{\i-1}
						\draw[edge] (u\j\a) to (u\i\a) {};
						\ifthenelse{\i<5}{
							\pgfmathtruncatemacro{\k}{\i+1}
							\draw[edge] (u\k\a) to (u\i\a) {};
						}{}
					}
				\end{scope}
			}
			\draw[edge] (v10) to (v61);
			\draw[edge] (v11) to (v60);
			\node[] (H) at (-4, -2.15) {$H:$};
		\end{tikzpicture}
		\caption{An illustration of an A-configuration $(D,X,H,f)$.}
		\label{fig:a_configuration}
	\end{center}
\end{figure}
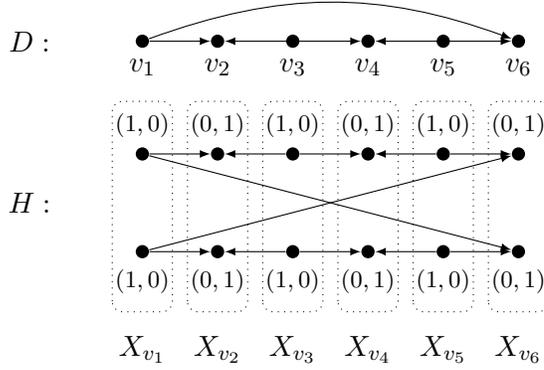

\item There are two disjoint constructible configurations, say  $(D^1,X^1,H^1,f^1)$ and $(D^2,X^2,H^2,f^2)$, such that $|D^i|<|D|$ for $i\in \{1,2\}$, $D$ is obtained from the disjoint graphs $D^1$ and $D^2$ by identifying a vertex $v^1 \in V(D^1)$ and a vertex $v^2 \in V(D^2)$ to a new vertex $v^\star$, $H$ is obtained from the disjoint digraphs $H^1$ and $H^2$ by choosing a bijection $\pi$ from $X_{v^1}$ to $X_{v^2}$ and identifying each vertex $x \in X_{v_1}$ with $\pi(x)$ to a vertex $x^\star$, and $f$ is defined as
$$f(y)=\begin{cases}
f^1(y) & \text{if } y \in V(H^1) \setminus X_{v^1},\\
f^2(y) & \text{if } y \in V(H^2) \setminus X_{v^2},\\
f^1(x) + f^2(\pi(x)) &\text{if } y \text{ is obtained from the identification of } \\ & x \in X_{v^1} \text{ with } \pi(x) \in X_{v^2}.
\end{cases}$$
In this case, we say that $(D,X,H,f)$ is obtained from $(D^1,X^1,H^1,f^1)$ and $(D^2,X^2,H^2,f^2)$ by \DF{merging} $v^1$ and $v^2$ to $v^\star$.
\end{enumerate}

By a \upshape{C}-\DF{configuration} we mean either an odd or an even \upshape{C}-configuration. An easy inductive argument imply that every constructible configuration is degree-feasible. The next result characterizes uncolorable degree-feasible configurations whose underlying graph is connected.

\begin{theorem} \label{main_theorem}
Let $D$ be a connected  digraph, let $(X,H)$ be a cover of $D$, and let $f:V(H) \to \mathbb{N}_0^2$ be a function. Then, $(D,X,H,f)$ is an uncolorable degree-feasible configuration if and only if $(D,X,H,f)$ is constructible.
\end{theorem}

The easy direction of the theorem, stating that constructible configurations are uncolorable, is proved in Section~\ref{subsec:struct_constructible}. The challenging direction is proved in Section~\ref{sectsub:proofs}. Our proof is inspired by the one of~\cite{GoncalvesPR2024}. The main idea is to consider a minimum counter-example $\cK = (D,X,H,f)$ and to boil down its possible structure, until we finally manage to explicitly give a transversal of $(X,H)$.

\subsection{Structure of the constructible configurations}
\label{subsec:struct_constructible}

The two next propositions will be useful in the proof of Theorem~\ref{main_theorem} when we want to justify that some configuration $\cK = (D,X,H,f)$ with $\suc(f) \neq \ems$ is indeed constructible. They both can be proved easily by induction, we omit the proofs.

 \begin{proposition}
 	\label{prop_add_0_0_arc}
 	Let $\cK=(D,X,H,f)$ be a constructible configuration, $uv$ be an arc of $D$, and $x_u \in X_u$, $x_v \in X_v$ be such that $A_H(\{x_u\}, X_v) = A_H(X_u,\{x_v\}) =  \ems$ and $f(x_u) = f(x_v) = (0,0)$. Then $\cK'=(D,X,H\cup \{x_ux_v\},f)$ is constructible.
  \end{proposition}

   \begin{proposition}
  	\label{prop_add_0_0_vertices}
  	Let $\cK=(D,X,H,f)$ be a constructible configuration and $\cK'=(D,X',H',f')$ be a configuration such that for some $v\in V(D)$ and some vertex $x_v$ of $H'$ we have:
  	\begin{itemize}
  		\item $x_v$ has no neighbor in $H'$,
  		\item $H = H' - x_v$,
  		\item $X'_v = X_v \cup \{x_v\}$ and $X'_u = X_u$ for every $u \in V(D)\setminus \{v\}$, and
  		\item $f'(x_v) = (0,0)$ and $f'(x) = f(x)$ for every $x \in V(H)\setminus \{x_v\}$.
  	\end{itemize}
  	Then $\cK'$ is constructible.
  \end{proposition}

The next proposition proves the ``if''-direction of Theorem~\ref{main_theorem}.

\begin{proposition}
\label{prop_constructible}
Let $(D,X,H,f)$ be a constructible configuration. Then the following statements hold:
\begin{enumerate}[label={\rm (\alph*)}]
\item $f(X_v) = d_D(v)$ for all $v \in V(D)$.
\label{enum:prop_constructible:a}
\item $(D,X,H,f)$ is uncolorable.
\label{enum:prop_constructible:b}
\end{enumerate}
\end{proposition}

\begin{proof}
Statement~\ref{enum:prop_constructible:a} can be proved easily by induction on the number of blocks of $D$. We prove~\ref{enum:prop_constructible:b} by contradiction. So assume that there exists a configuration $\cK=(D,X,H,f)$ such that  $\cK$ is constructible and colorable, i.e., there is a transversal $T$ of $(X,H)$ such that $H[T]$ is strictly $f$-degenerate. Among all such configurations, we let $\cK$ be the one for which $|D|$ is minimum.

Recall that if $f(x)=(0,0)$ for some vertex $x \in V(H)$, then $H[\{x\}]$ is not strictly $f$-degenerate and, hence, $x$ cannot be contained in any strictly $f$-degenerate subdigraph of $H$. Thus, $T\subseteq \su(f)$.

First, assume that $(D,X,H,f)$ is an \upshape{M}-configuration. Then $D$ is a block and there exists a transversal $U$ of $(X,H)$ such that $(X,H)/U$ is a saturated cover of $D$. Furthermore, for $v\in V(D)$ and $x\in X_v$, we have $f(x)=d_D(v)$ if $x\in U$ and $f(x)=(0,0)$ otherwise. Consequently, $T = U$ and $d_{H[T]}(x)=d_D(v)=f(x)$ for all $x\in T$, a contradiction.

Next assume that $(D,X,H,f)$ is a \upshape{K}-configuration. Then $D=D^\pm(K_n)$ with $n\geq 2$ and there are integers $n_1,n_2,\ldots,n_p \in \nat$ with $p \geq 1$ such that $n_1 + n_2 + \cdots + n_p=n-1$. Moreover, there are $p$ disjoint transversals $U_1, U_2, \ldots, U_p$ of $(X,H)$ such that,  for $i\in [1,p]$, $(X,H)/U_i$ is a saturated cover of $D$. Furthermore, for $v\in V(H)$ and $x\in X_v$, we have $f(x)=(n_i,n_i)$ if $x\in U_i$ for $i\in [1,p]$ and $f(x)=(0,0)$ otherwise. Then $T$ is a subset of $U_1 \cup U_2 \cup \cdots \cup U_p$. Since $|T|=n$, there exists $i\in [1,p]$, such that $m=|T\cap U_i|\geq n_i+1$. Then $H'=H[T\cap U_i]$ is a copy of $D^\pm(K_m)$ with $m\geq n_i+1$ implying that $d_{H'}(x)=(m-1,m-1) \geq (n_i,n_i) =f(x)$ for every $x\in V(H')$, a contradiction. If $(D,X,H,f)$ is a \upshape{C}-configuration or an \upshape{A}-configuration, the argument is similar.

Hence, consider the case that $\cK=(D,X,H,f)$ is obtained from two constructible configurations, say $\cK^1=(D^1,X^1,H^1,f^1)$ and $\cK^2=(D^2,X^2,H^2,f^2)$, by merging $v^1 \in V(D^1)$ and $v^2 \in V(D^2)$ to a new vertex $v^\star$. For the sake of clarity, we assume that $v^1=v^2=v^\star$ and $X_{v^\star}^1=X_{v^\star}^2$ (i.e., $\pi=id$). By minimality of $|D|$, we conclude that $\cK^i$ is uncolorable for $i \in \{1,2\}$. 

Let $T$ be a transversal of $(X,H)$ such that $H[T]$ is strictly $f$-degenerate. Let $T^i=T\cap V(H^i)$ ($i\in \{1,2\}$) and let $x^\star$ be the unique vertex from $X_{v^\star}\cap T$. Then $T^1\cap T^2=\{x^\star\}$, $f(x^\star)=f^1(x^\star) + f^2(x^\star)$, and $f(x)=f^i(x)$ for all $x \in T^i \sm \{x^\star\}$ ($i \in \{1,2\}$). Since $\cK^i$ is uncolorable, the subgraph $H[T^i]$ is not strictly $f^i$-degenerate, implying that there is a nonempty subdigraph $\tilde{H^i}$ of $H[T^i]$ such that $d_{\tilde{H^i}}(x) \geq f^i(x)$ for all $x \in V(\tilde{H^i})$ ($i\in \{1,2\}$). If $x^\star$ does not belong to $\tilde{H^i}$, then $\tilde{H^i}$ is a subgraph of $H[T]-x^\star$ and so
$\tilde{H^i}$ is strictly $f^i$-degenerate as $H[T]$ is strictly $f$-degenerate, a contradiction. Hence $x^\star$ belongs to both $\tilde{H^1}$ and $\tilde{H^2}$. Let $\tilde{H}$ be the subdigraph of $H[T]$ induced by $V(\tilde{H^1}) \cup V(\tilde{H^2})$. Let $x$ be an arbitrary vertex of $\tilde{H}$. If $x\not=x^\star$, then $x$ belongs to $\tilde{H^i}-x^\star$ for some $i\in \{1,2\}$, and so $d_{\tilde{H}}(x)=d_{\tilde{H^i}}(x) \geq f^i(x) = f(x)$. Furthermore, we have
$$d_{\tilde{H}}(x^\star)=d_{\tilde{H^1}}(x^\star) + d_{\tilde{H^2}}(x^\star) \geq f^1(x^\star) + f^2(x^\star) = f(x^\star).$$
Hence, $\tilde{H} \subseteq H[T]$ is not strictly $f$-degenerate and so $H[T]$ is not strictly $f$-degenerate as well, a contradiction.
\end{proof}

\subsection{Proof of Theorem~\ref{main_theorem}}
\label{sectsub:proofs}

To complete the proof of Theorem~\ref{main_theorem}, it remains to show that each uncolorable degree-feasible configuration is constructible. To this end, we will apply the following reduction method many times.

\begin{proposition}[General reduction]\label{prop_reduction_general}
Let $\cK=(D,X,H,f)$ be a configuration. Let $T$ be a partial transversal of $(X,H)$ such that $H[T]$ is strictly $f$-degenerate and $D' = D - \dom(T:D)$ is connected. Let $\cK'=(D',X',H',f')$ be the configuration with $(X',H')=(X,H)/D'$, and let
$$f'^+(x)=\max \Big\{0,f^+(x)-\sum_{y \in T} a_H(x,y) \Big\} \mbox{ and } f'^-(x)=\max  \Big\{0,f^-(x)-\sum_{y \in T} a_H(y,x) \Big\}$$
for all $x\in V(H')$.
Then the following statements hold:
\begin{enumerate}[label={\rm (\alph*)}]
\item If $\cK$ is degree-feasible, then so is $\cK'$.
\label{enum:prop_reduction_general:a}
\item If $\cK$ is uncolorable, then so is $\cK'$.
\label{enum:prop_reduction_general:b}
\end{enumerate}
In the following, we write $\cK'=\cK/T$. If further $T=\{x_v\}$ with $x_v\in X_v$, we write $\cK'=\cK/(v,x_v)$.
\end{proposition}

\begin{proof}
For the proof of~\ref{enum:prop_reduction_general:a}, assume that $\cK$ is degree-feasible. Let $u$ be an arbitrary vertex of $D'$. Using Proposition~\ref{prop:cover}, we have
\begin{eqnarray*}
f'^+(X_u) &=& \sum_{x\in X_u} f'^+(x)\geq \sum_{x\in X_u} \Big(f^+(x) - \sum_{y \in T} a_H(x,y)\Big)\\
          &=& \sum_{x\in X_u} f^+(x)-\sum_{y \in T}\sum_{x\in X_u} a_H(x,y)\\
          &\geq& d^+_D(u)-\left|A_D\big(\{u\},\dom(T:D)\big)\right|=d^+_{D'}(u).
\end{eqnarray*}
Similarly we have $f'^-(X_u) \geq d^-_{D'}(u)$.
Consequently, $\cK'$ is degree-feasible, which proves~\ref{enum:prop_reduction_general:a}. 

For the proof of~\ref{enum:prop_reduction_general:b}, assume that $\cK'$ is colorable, so there is a transversal $T'$ of $(X',H')$ such that $H'[T']$ is strictly $f'$-degenerate. Then $T_0=T \cup T'$ is a transversal of $(X,H)$ and we claim that $H[T_0]$ is strictly $f$-degenerate, which implies that $\cK$ is colorable as well. To prove the claim, let $\tilde{H}$ be an arbitrary nonempty subdigraph of $H[T_0]$, and let $L=\tilde{H}-T$. If $L=\ems$, then $\tilde{H}$ is a nonempty subdigraph of $H[T]$. Since $H[T]$ is strictly $f$-degenerate, we are done. Otherwise, $L$ is a nonempty subdigraph of $H'[T']=H[T]$. Since $H'[T']$ is strictly $f'$-degenerate, there is a vertex $x$ in $L$ such that  $d^+_L(x)<f'^+(x)$ (or $d^-_L(x)<f'^-(x)$). This implies that $f'^+(x)=f^+(x)- \sum_{y \in T} a_H(x,y)$ (or $f'^-(x)=f^-(x)- \sum_{y \in T} a_H(y,x)$), which gives
$d^+_{\tilde{H}}(x)\leq d^+_L(x)+ \sum_{y \in T} a_H(x,y)<f^+(x)$ (or $d^-_{\tilde{H}}(x)\leq d^-_L(x)+ \sum_{y \in T} a_H(y,x)<f^-(x)$), concluding the proof of~\ref{enum:prop_reduction_general:b}.
\end{proof}

The following is the direct application of Proposition~\ref{prop_reduction_general} for $T=\{x_v\}$ with $x_v\in X_v$.
\begin{proposition}[Simple reduction]\label{prop_reduction}
	Let $\cK=(D,X,H,f)$ be a configuration, let $v \in V(D)$ be a non-separating vertex and $x_v \in X_v$ be such that $f(x_v)\not=(0,0)$. Then the following statements hold:
	\begin{enumerate}[label={\rm (\alph*)}]
		\item If $\cK$ is degree-feasible, then so is $\cK/(v,x_v)$.
		\label{enum:prop_reduction:a}
		\item If $\cK$ is uncolorable, then so is $\cK/(v,x_v)$.
		\label{enum:prop_reduction:b}
	\end{enumerate}
\end{proposition}

Using the reduction method, we obtain the following useful properties of uncolorable degree-feasible configurations.

\begin{proposition}\label{prop_uncolorable-facts}
Let $D$ be a connected digraph and let $(D,X,H,f)$ be an uncolorable degree-feasible configuration. Then the following statements hold:
\begin{enumerate}[label={\rm (\alph*)}]
\item $f(X_u)=d_D(u)$ for all $u \in V(D)$.
\label{enum:prop_uncolorable-facts:a}

\item 
\label{enum:prop_uncolorable-facts:b}
Suppose that $v\in V(D)$ is non-separating, $x_v \in X_v$ satisfies $f(x_v)\not=(0,0)$, and $u\in V(D-v)$. Then the following statements hold:
    \begin{enumerate}[label=\arabic*.]
      \item $f^+(x)\geq a_H(x,x_v)$ for all $x\in X_u$ and $\sum_{x\in X_u} a_H(x,x_v)=a_D(u,v)$.
      \item $f^-(x)\geq a_H(x_v,x)$ for all $x\in X_u$ and $\sum_{x\in X_u} a_H(x_v,x)=a_D(v,u)$.
    \end{enumerate}
\item If $|D| \geq 2$ and if $u$ is an arbitrary vertex of $D$, then there is a partial transversal $T$ of $(X,H)$ such that $\dom(T:H)=V(D-u)$ and $H[T]$ is strictly $f$-degenerate. Furthermore, for every such transversal $T$ and every vertex $x\in X_u$, we have $f(x)=d_{H[T\cup \{x\}]}(x)$.
\label{enum:prop_uncolorable-facts:c}
\end{enumerate}
\end{proposition}
\begin{proof}
We prove~\ref{enum:prop_uncolorable-facts:a} and~\ref{enum:prop_uncolorable-facts:b} simultaneously by induction on the order of $D$. If $|D|=1$, then $V(H)=\suc(f)$, and both statements are evident. Suppose that $|D|\geq 2$, and let $u$ be an arbitrary vertex of $D$.
Then $D$ has a nonseparating vertex $v$ distinct from $u$. Since $f(X_v)\geq d_D(v)\not=(0,0)$, there is a color $x_v\in X_v$ such that $f(x_v)\not=(0,0)$. From Proposition~\ref{prop_reduction} it follows that $(D',X',H',f')=(D,X,H,f)/(v,x_v)$ is an uncolorable degree-feasible configuration; note that $D'=D-v$, $H'=H-X_v$ and  $X'_w=X_w$ for all vertices $w$ of $D'$; furthermore, we have
$$f'^+(x)=\max\{0,f^+(x)-a_H(x,x_v)\} \mbox{ and } f'^-(x)=\max\{0,f^-(x)-a_H(x_v,x)\}$$
for all $x\in V(H')$.
By the induction hypothesis, we obtain that
$$f'(X_u)=f'(X'_u)=d_{D'}(u).$$
Using Proposition~\ref{prop:cover}, we then obtain that
\begin{eqnarray*}
d^+_{D'}(u)=f'^+(X_u) &=& \sum_{x\in X_u} f'^+(x)\geq \sum_{x\in X_u}(f^+(x)-a_H(x,x_v))\\
          &=& \sum_{x\in X_u} f^+(x)-\sum_{x\in X_u}a_H(x,x_v)\\
          &\geq& f^+(X_u)-a_D(u,v)\\
          &\geq& d^+_D(u)-a_D(u,v)=d^+_{D'}(u)
\end{eqnarray*}
This implies that  $f^+(X_u) = d^+_{D'}(u) + a_D(u,v) = d^+_D(u)$, $f^+(x)\geq a_H(x,x_v)$ for all $x\in X_u$, and $\sum_{x\in X_u} a_H(x,x_v)=a_D(u,v)$. The same arguments imply  $f^-(X_u)=d^-_D(u)$, $f^-(x)\geq a_H(x_v,x)$ for all $x\in X_u$ and $\sum_{x\in X_u} a_H(x_v,x)=a_D(v,u)$.
This proves~\ref{enum:prop_uncolorable-facts:a} and~\ref{enum:prop_uncolorable-facts:b}.

For the proof of~\ref{enum:prop_uncolorable-facts:c}, suppose that $|D| \geq 2$ and $u \in V(D)$. Let $D'=D-u$, let $X'$ be the restriction of $X$ to $V(D')$, and let $H'=H-X_u$. Then $(D',X',H',f)$ is a degree-feasible configuration. Since $D$ is connected, each component of $D'$ contains a vertex $v\in N^+_D(u) \cup N^-_D(v)$, which thus satisfies
$$f^+(X_v) > d^+_{D'}(v) \mbox{ or } f^-(X_v) > d^-_{D'}(v).$$
Then it follows from~\ref{enum:prop_uncolorable-facts:a} applied to each component of $D'$ that $(X',H')$ has a transversal $T$ such that $H'[T]$ is strictly $f$-degenerate. Then $T$ is a partial transversal  of $(X,H)$ such that $\dom(T:H)=V(D')$ and $H[T]=H'[T]$ is strictly $f$-degenerate. Now let $T$ be such a partial transversal and let $x\in X_u$ be an arbitrary vertex. Since $(D,X,H,f)$ is uncolorable, $H[T \cup \{x\}]$ is not strictly $f$-degenerate. Since $H[T]$ is strictly $f$-degenerate, we then conclude that $d_{H[T \cup \{x\}]}(x) \geq f(x)$. Using \ref{enum:prop_uncolorable-facts:a} and \ref{desc:D2}, we then obtain that
$$d_D(u)=f(X_u)=\sum_{x\in X_u} f(x)\leq \sum_{x\in X_u}d_{H[T\cup\{x\}]}(x)\leq d_D(u).$$
This obviously implies that $d_{H[T \cup \{x\}]}(x) = f(x)$ for all $x\in X_u$, which proves~\ref{enum:prop_uncolorable-facts:c}.
\end{proof}

We will use the following decomposition of 2-connected undirected graphs, similar to the well-known ear decompositions.

\begin{lemma}[{\cite[Lemma~26]{GoncalvesPR2024}}]
	\label{lemma:decomposition}
	If $G$ is a 2-connected graph, then at least one of the following holds:
	\begin{itemize}
		\item $G$ is a cycle on at least three vertices,
		\item $G$ contains a vertex $v$ such that $G-v$ is 2-connected, or
		\item $G$ contains a path $P$ on at least two vertices such that
		$G - V(P)$ is 2-connected and each vertex of $P$ has degree exactly two in $G$.
	\end{itemize}
\end{lemma}

Next theorem proves the ``only if''-direction of Theorem~\ref{main_theorem}.

\begin{theorem} \label{theorem:main2}
If $\cK=(D,X,H,f)$ is an uncolorable degree-feasible configuration for which $D$ is connected, then $\cK$ is constructible.
\end{theorem}

\begin{proof}
The proof is by contradiction. Let $\cK=(D,X,H,f)$ be a minimal counter-example, that is, $D$ is a connected digraph such that
\begin{enumerate}[label={\rm (\Alph*)}]
\item $\cK$ is an uncolorable degree-feasible configuration,
\label{enum:proof:thm2:A}
\item $\cK$ is not constructible,
\label{enum:proof:thm2:B}
\item $|D|$ is minimum subject to~\ref{enum:proof:thm2:A} and~\ref{enum:proof:thm2:B}, and
\label{enum:proof:thm2:C}
\item $|H| + |A(H)|$ is minimum subject to~\ref{enum:proof:thm2:A},~\ref{enum:proof:thm2:B}, and~\ref{enum:proof:thm2:C}.
\label{enum:proof:thm2:D}
\end{enumerate}
By Proposition~\ref{prop_uncolorable-facts}\ref{enum:prop_uncolorable-facts:a} we have
\begin{align}\label{eq_degree-sum}
f(X_v) = d_D(v) \mbox{ for all } v \in V(D).
\end{align}
Clearly, $|D| \geq 2$, as for $|D|=1$ we have $V(D)=\{v\}$ and $f(x)=(0,0)$ for all $x \in X_v$ implying that $(D,X,H,f)$ is an \upshape{M}-configuration and hence constructible, a contradiction to~\ref{enum:proof:thm2:B}. We will reach a contradiction via a sequence of claims. We first prove that $D$ is a block.

\begin{claim} \label{claim:block}
$D$ is a block.
\end{claim}
\begin{proofclaim}
	Assume for a contradiction that $D$ is the union of two connected subdigraphs $D^1$ and $D^2$ such that $V(D^1)\cap V(D^2)=\{v^\star\}$ and $|D^i| < |D|$ for $i \in \{1,2\}$. For $i\in\{1,2\}$, let $(X^i,H^i)=(X,H)/D^i$. We now define a vertex function $f^i$ of $H^i$ as follows. By Proposition~\ref{prop_uncolorable-facts}\ref{enum:prop_uncolorable-facts:c}, $(X,H)$ has a partial transversal $T$ such that $\dom(T:D)=V(D-v^\star)$ and $H[T]$ is strictly $f$-degenerate. Let $T_1=T \cap V(H^1)$ and let $T_2= T \cap V(H^2)$. Then $H[T]$ is the disjoint union of $H[T_1]$ and $H[T_2]$; note that $H[T_i]=H^i[T_i]$. Then, using Proposition~\ref{prop_uncolorable-facts}\ref{enum:prop_uncolorable-facts:c}, we obtain that
$$f(x) = d_{H[T \cup \{x\}]}(x)=d_{H[T_1 \cup \{x\}]}(x) + d_{H[T_2 \cup \{x\}]}(x)$$
for all $x \in X_{v^\star}$. We set $f^i(x)=d_{H[T_i \cup \{x\}]}(x)$ for $i \in \{1,2\}$ and $x \in X_{v^\star}$. For a vertex $v$ of $D^i-v^\star$, let $f^i(x)=f(x)$ for all $x\in X_v$. Clearly, $\cK^i=(D^i,X^i,H^i,f^i)$ is a configuration for $i\in \{1,2\}$. We show that, for $i\in \{1,2\}$, $\cK^i$ is both uncolorable and degree-feasible.

We first show that $\cK^1$ and $\cK^2$ are uncolorable. Assume not, so by symmetry there is a transversal $T^1$ of $(X^1,H^1)$ such that $H^1[T^1]$ is strictly $f^1$-degenerate. Clearly, $T^1 \cup T_2$ is a transversal of $(X,H)$, and we claim that $H[T^1\cup T_2]$ is strictly $f$-degenerate. If this is not the case, then there is a nonempty subdigraph $\tilde{H}$ of $H[T^1\cup T_2]$ with $d_{\tilde{H}}(x) \geq f(x)$ for all $x \in V(\tilde{H})$. Since $H[T_2]$ is strictly $f$-degenerate, $\tilde{H}$ contains vertices of $T^1$. Since $H[T^1]$ is strictly $f^1$-degenerate, there is a vertex $y\in V(\tilde{H})\cap T^1$ such that
$d^+_{\tilde{H}-T_2}(y) < {f^1}^+(y)$ or $d^-_{\tilde{H}-T_2}(y) < {f^1}^-(y)$.
If $y \not \in X_{v^\star}$, then we obtain that
$$f^+(y)\leq d^+_{\tilde{H}}(y)= d^+_{\tilde{H}-T_2}(y) < {f^1}^+(y) = f^+(y)$$
or
$$f^-(y)\leq d^-_{\tilde{H}}(y)= d^-_{\tilde{H}-T_2}(y) < {f^1}^-(y) = f^-(y),$$
a contradiction. If $y \in X_{v^\star}$, then $f^2(y)=d_{H[T_2 \cup \{y\}]}(y)\geq d_{\tilde{H}-T^1}(y)$ and we obtain that
$$f^+(y)\leq d^+_{\tilde{H}}(y)= d^+_{\tilde{H}-T_2}(y) + d^+_{\tilde{H}-T^1}(y) < {f^1}^+(y) + {f^2}^+(y) = f^+(y)$$
or
$$f^-(y)\leq d^-_{\tilde{H}}(y)= d^-_{\tilde{H}-T_2}(y) + d^-_{\tilde{H}-T^1}(y) < {f^1}^-(y) + {f^2}^-(y) = f^-(y),$$
which is a contradiction again. Hence, $H[T^1\cup T_2]$ is strictly $f$-degenerate and so $\cK$ is colorable, a contradiction to~\ref{enum:proof:thm2:A}. This shows that $\cK^1$ and $\cK^2$ are uncolorable.

Let us now show that $\cK^i$ is degree-feasible for $i\in \{1,2\}$. By \eqref{eq_degree-sum} and the definition of $f^1$, we obtain that $f^1(X_v)=f(X_v)=d_D(v)=d_{D^1}(v)$
for all $v \in V(D^1-v^\star)$. Moreover, we have
\begin{align}\label{eq_degree_feasible}
d_{D^1}(v^\star) + d_{D^2}(v^\star) = d_D(v^\star)  = f(X_{v^\star}) =f^1(X_{v^\star}) + f^2(X_{v^\star}).
\end{align}
For every vertex $y\in T$ there is a unique vertex $v_y$ in $D$ such that $y\in X_{v_y}$. Using Proposition~\ref{prop:cover}, we then obtain for $i\in \{1,2\}$ that

\begin{eqnarray*}
{f^i}^+(X_{v^\star})&=&\sum_{x\in X_{v^\star}}d^+_{H[T_i\cup\{x\}]}(x)=\sum_{x\in X_{v^\star}}\sum_{y\in T_i}a_H(x,y)\\
&=& \sum_{y\in T_i}\sum_{x\in X_{v^\star}}a_H(x,y) \leq \sum_{y\in T_i}a_D(v^\star,v_y)=d^+_{D_i}(v^\star).
\end{eqnarray*}
Similarly, we obtain than ${f^i}^-(X_{v^\star})\leq d^-_{D_i}(v^\star)$ for $i\in \{1,2\}$.
By \eqref{eq_degree_feasible}, this implies that
$f^i(X_{v^\star})= d_{D^i}(v^\star)$
for $i\in \{1,2\}$. Consequently, $\cK^i$ is degree-feasible for $i\in \{1,2\}$.

Since $\cK=(D,X,H,f)$ is a minimal counter-example and $|D^i|<|D|$ for $i\in \{1,2\}$, we conclude that $\cK^i=(D^i,X^i,H^i,f^i)$ is a constructible configuration, and so $\cK$ is obtained from the constructible configurations $\cK^1$ and $\cK^2$ by merging two vertices to $v^\star$. Hence, $\cK$ is a constructible configuration, a contradiction to~\ref{enum:proof:thm2:B}.
\end{proofclaim}

\begin{claim}
	\label{claim:0_0_isolated}
	If there exists $x\in V(H)$ with $f(x) = (0,0)$, then $x$ has no neighbor in $H$.
\end{claim}

\begin{proofclaim}
	Assume for a contradiction that there exists $x\in V(H)$ with $f(x) = (0,0)$ and that $x$ has a neighbor $y$ in $H$. Since $D$ is a block (by Claim~\ref{claim:block}), we have $f(y) = (0,0)$ (by Proposition~\ref{prop_uncolorable-facts}\ref{enum:prop_uncolorable-facts:b}). Let $H' = H - \{xy,yx\}$ and $\cK' = (D, X, H', f)$. Clearly, $\cK'$ is a configuration. Since $\cK$ is uncolorable and degree-feasible (by 
\ref{enum:proof:thm2:A}), we obtain that $\cK'$ is uncolorable and degree-feasible, too. Since $|H'|+|A(H')|<|H|+|A(H)|$, it then follows from
\ref{enum:proof:thm2:C} and
\ref{enum:proof:thm2:D} that $\cK'$ is constructible. This implies that $\cK$ is constructible (by Proposition~\ref{prop_add_0_0_arc}), a contradiction to 
\ref{enum:proof:thm2:B}. 
\end{proofclaim}

\begin{claim}
	\label{claim:no_0_0}
	For every vertex $x\in V(H)$, $f(x) \neq (0,0)$.
\end{claim}
\begin{proofclaim}
	Assume for a contradiction that there exists $x\in V(H)$ with $f(x) = (0,0)$. Then $x$ has no neighbor in $H$ (by Claim~\ref{claim:0_0_isolated}). Let $(X',H')=(X,H)/V(H-x)$ and let $f'$ be the restriction of $f$ to $V(H')=V(H-x)$. Then $\cK'=(D,X',H',f')$ is an uncolorable degree-feasible configuration (by 
\ref{enum:proof:thm2:A}). Since $|H'|+|A(H')|<|H|+|A(H)|$, it then follows from
\ref{enum:proof:thm2:C} and
\ref{enum:proof:thm2:D} that $\cK'$ is constructible. This implies that $\cK$ is constructible (by Proposition~\ref{prop_add_0_0_vertices}), a contradiction to 
\ref{enum:proof:thm2:B}. 
\end{proofclaim}

\begin{claim} \label{claim:uniform}
	$(X,H)$ is a $r$-uniform cover of $D$ for some $r\geq 2$.
\end{claim}
\begin{proofclaim}
	First note that $(X,H)$ is not $1$-uniform, for otherwise $\cK$ is an \upshape{M}-configuration (by Claim~\ref{claim:block} and \eqref{eq_degree-sum}), a contradiction.
	
	Assume for a contradiction that $(X,H)$ is not $r$-uniform for any $r\geq 2$. Since $D$ is connected, there exist two adjacent 
vertices $u,v \in V(D)$ such that $|X_u| < |X_v|$. Recall that, in $H$, both the arcs from $X_u$ to $X_v$  and the arcs from $X_v$ to $X_u$ form a matching. Hence, there exist $x_v, y_v \in X_v$ satisfying
		\[ \sum_{x \in X_u} a_H(x,x_v) = 0 \mbox{~~~and~~~} \sum_{x \in X_u} a_H(y_v,x) = 0. \]
Note that $D$ is a block (by Claim~\ref{claim:block}) and we have $f(x_v) \neq (0,0)$ as well as $f(y_v) \neq (0,0)$ (by Claim~\ref{claim:no_0_0}). It then follows from 
\ref{enum:proof:thm2:A} and Proposition~\ref{prop_uncolorable-facts}\ref{enum:prop_uncolorable-facts:b} that $a_D(u,v) = a_D(v,u)=0$, a contradiction.
\end{proofclaim}

\begin{claim} \label{claim:perfect_matching}
Let $uv\in A(D)$. Then the following statements hold:
\begin{enumerate}[label={\rm (\alph*)}]
\item $A_H(X_u,X_v)$ is a perfect matching of $H[X_u \cup X_v]$.
\label{enum:claim:perfect_matching:a}
\item $f^+(x_u)\geq 1$ for every $x_u\in X_u$ and $f^-(x_v)\geq 1$ for every $x_v\in X_v$.
\label{enum:claim:perfect_matching:b}
\item $d_D^+(u)=f^+(X_u)\geq |X_u|=r$ and $d_D^-(v)=f^-(X_v)\geq |X_v|= r$.
\label{enum:claim:perfect_matching:c}
\end{enumerate}
\end{claim}
\begin{proofclaim}
By Claim~\ref{claim:uniform}, we have $|X_u|=|X_v|=r$. Since $uv\in A(D)$, $A_H(X_u,X_v)$ is a matching of $H[X_u \cup X_v]$ (by \ref{desc:D2}). If this matching is not perfect, then there exists $x_u \in X_u$ satisfying
	\[ \sum_{x \in X_v} a_H(x_u,x) = 0. \]
	By Claim~\ref{claim:no_0_0}, we have $f(x) \neq (0,0)$ for every $x\in V(H)$. Thus Proposition~\ref{prop_uncolorable-facts}\ref{enum:prop_uncolorable-facts:b} implies that $a_D(u,v) =0$, a contradiction. This proves \ref{enum:claim:perfect_matching:a}. Since $D$ is a block (by Claim~\ref{claim:block}), statement \ref{enum:claim:perfect_matching:b} follows from 
\ref{enum:claim:perfect_matching:a} and Proposition~\ref{prop_uncolorable-facts}\ref{enum:prop_uncolorable-facts:b}. Statement \ref{enum:claim:perfect_matching:c} follows from \ref{enum:claim:perfect_matching:b} and \eqref{eq_degree-sum}. 
\end{proofclaim}

\begin{claim}
	\label{claim:2connected}
	$D$ is 2-connected.
\end{claim}
\begin{proofclaim}
	Recall that $|D| \geq 2$ and $D$ is a block (by Claim~\ref{claim:block}). Hence, if $D$ is not 2-connected, then $|D| = 2$ and $D$ contains an arc $uv$. Using Claims~\ref{claim:uniform} and \ref{claim:perfect_matching}\ref{enum:claim:perfect_matching:c}, we then obtain that $1=d_D^+(u)=f^+(X_u) \geq r \geq 2$, which is impossible.
\end{proofclaim}

\begin{claim}
\label{claim:reduction}
Let $v\in V(D)$ be an arbitrary vertex and let $x_v\in V(H)$ be an arbitrary vertex. Then $\cK'=\cK/(v,x_v)$ is a constructible configuration. Furthermore, if $\cK'=(D',X',H',f')$ then the following statements hold:	
\begin{enumerate}[label={\rm (\alph*)}]
\item $D'=D-v$ and $(X',H')=(X,H)/D'$.
\label{enum:claim:reduction:a}
\item $f'^+(x)=f^+(x)-a_H(x,x_v)\geq 0$ for all $x\in V(H')$.
\label{enum:claim:reduction:b}
\item $f'^-(x)=f^-(x)-a_H(x_v,x)\geq 0$ for all $x\in V(H')$. 
\label{enum:claim:reduction:c}
\end{enumerate}
In what follows we say that $f'$ is the \DF{vertex function} of $\cK'$.
\end{claim}
\begin{proofclaim}
Note that $v$ is a nonseparating vertex of $D$ (by Claim~\ref{claim:2connected}) and $f(x_v)\not=(0,0)$ (by Claim~\ref{claim:no_0_0}). Then, it follows from
Proposition~\ref{prop_reduction} and \ref{enum:proof:thm2:A} that $\cK'$ is an uncolorable and degree-feasible configuration. By \ref{enum:proof:thm2:B} and \ref{enum:proof:thm2:C} this implies that $\cK'$ is a constructible configuration. Statement \ref{enum:claim:reduction:a} follows from the definition of $\cK'$. Statements (b) and (c) are consequences of the definition of $f'$ and Proposition~\ref{prop_uncolorable-facts}\ref{enum:prop_uncolorable-facts:b}.
\end{proofclaim}

In what follows, let $G$ denote the underlying graph of $D$. By Claim~\ref{claim:2connected}, $G$ is 2-connected. First, we investigate a vertex of degree two in $G$. Note that $|G|\geq 3$ and so $\de(G)\geq 2$. 

\begin{claim}
\label{claim:vertexG}
Let $u$ be a vertex of $G$ such that $d_G(u)=2$, and let $N_G(u)=\{v,v'\}$. Then the following statements hold:
\begin{enumerate}[label={\rm (\alph*)}]
\item $d_D^\pm(u)\in \{0,2\}$ and $r=2$.
\label{enum:claim:vertexG:a}
\item \label{enum:claim:vertexG:b}
If $D[u,v]$ is digon, then $D[u,v']$ is a digon, $f(x)=(1,1)$ for every $x\in X_u$, and $H[X_u \cup X_v]$ is a bidirected graph. 
\item If $D[u,v]$ is not a digon, then $u$ is either a source of $D$ and $f(x)=(1,0)$ for every $x\in X_u$, or $u$ is a sink of $D$ and $f(x)=(0,1)$ for every $x\in X_u$. 
\label{enum:claim:vertexG:c}
\end{enumerate}
\end{claim}
\begin{proofclaim}
Since $d_G(u)=2$, we have $d_D^\pm(u)\in \{0,1,2\}$. Using Claims \ref{claim:uniform} and \ref{claim:perfect_matching}, we obtain 
\ref{enum:claim:vertexG:a}. For the proof of \ref{enum:claim:vertexG:b}, assume that $D[u,v]$ is a digon. Based on \ref{enum:claim:vertexG:a} and Claim~\ref{claim:perfect_matching}, we obtain that $D[u,v']$ is a digon, too, and $f(x)=(1,1)$ for every $x\in X_u$. Suppose, to the contrary, that $H_1=H[X_u\cup X_v]$ is not a bidirected graph. Since $|X_u|=|X_v|=r=2$ it then follows from Claim~\ref{claim:perfect_matching}, that $H_1$ is a directed cycle on four vertices, say $C=(x_v,x_u,y_v,y_u,x_v)$ with $X_u=\{x_u,y_u\}$ and $X_v=\{x_v,y_v\}$. Now let $\cK'=(D',X',H',f')=\cK/(v,x_v)$. Using Claim~\ref{claim:reduction}, we obtain that $\cK'$ is a constructible configuration. Furthermore, we obtain that $f'(x_u)=(1,0)$, since $f(x_u)=(1,1)$. Then $B=D[u,v']$ is a block of $D'=D-v$. Clearly, $B$ is the only block of $D'$ containing $u$ and $B=D^\pm(K_2)$ (by \ref{enum:claim:vertexG:b}). Since $\cK'$ is constructible and $r=2$, this implies that $f'(x_u)\in \{(0,0),(1,1)\}$, a contradiction. This proves \ref{enum:claim:vertexG:b}. Using \ref{enum:claim:vertexG:a} and Claim~\ref{claim:perfect_matching}, we obtain 
\ref{enum:claim:vertexG:c}.
\end{proofclaim}

\begin{claim}
\label{claim:Gnotcycle}
The underlying graph $G$ of $D$ is not a cycle.
\end{claim}
\begin{proofclaim}
Suppose this is false, that is, $G$ is a cycle. By Claim~\ref{claim:vertexG}, this implies that $(X,H)$ is a $2$-uniform cover of $D$. First, assume that $D$ contains a digon. Then it follows from Claim~\ref{claim:vertexG} that $D$ is a bidirected cycle, say $D=D^\pm(C_n)$ with $n\geq 3$, $f(x)=(1,1)$ for every $x\in V(H)$, and $H$ is either a bidirected cycle of length $2n$, or the disjoint union of two bidirected cycles both of length $n$. In the former case, as $\cK$ is not constructible, $n$ must be odd. In the latter case, $n$ must be even. In both cases, it is straightforward that $(X,H)$ admits a  transversal $T$ which is an independent set of $H$. As $f(x)=(1,1)$ for every $x \in T$, we conclude that $\cK$ is colorable, a contradiction. Now, assume that $D$ contains no digon. Then it follows from Claim~\ref{claim:vertexG} that $D$ is an antidirected cycle. Furthermore, we obtain that if $u$ is a source of $D$, then $f(x)=(1,0)$ for every $x\in X_u$, and if $u$ is a sink of $D$, then $f(x)=(0,1)$ for every $x\in X_u$. Using Claim~\ref{claim:perfect_matching}\ref{enum:claim:perfect_matching:a}, we also obtain that 
$H$ is either an antidirected cycle on $2|D|$ vertices or the disjoint union of two antidirected cycles. In the former case, $\cK$ is an \upshape{A}-configuration, a contradiction. In the latter case, we can find a transversal of $(X,H)$ which is an independent set of $H$, which again yields a contradiction.
\end{proofclaim}

\begin{claim}
\label{claim:Dremains}
$G$ contains a vertex $v$ such that $G-v$ is 2-connected.
\end{claim}
\begin{proofclaim}
Suppose this is false. Since $G$ is not a cycle (by Claim~\ref{claim:Gnotcycle}), it then follows from Lemma~\ref{lemma:decomposition} that 
$G$ contains a path $P=(v_1, v_2,\dots,v_\ell)$ on $\ell \geq 2$ vertices such that $G - V(P)$ is 2-connected and each vertex of $P$ has degree exactly two in $G$.

Let $v_0$ be the only neighbor of $v_1$ distinct from $v_2$  in $G$ and $v_{\ell+1}$ be the only neighbor of $v_\ell$ distinct from $v_{\ell - 1}$ in $G$. Note that $v_0$ is distinct from $v_{\ell+1}$ because $D$ is 2-connected. Since $G-V(P)$ is 2-connected, there is a vertex $u$ in $G$
distinct from $v_0, v_1, \dots,v_{\ell+1}$. Note that $u$ is, in $D$, not adjacent to any vertex of $V(P)$. Let $D'=D[\{v_0,v_1, \ldots, v_{\el+1}\}]$ and $(X',H')=(X,H)/D'$. Using Claim~\ref{claim:vertexG}, we obtain that $(X,H)$ is a 2-uniform cover of $D$, and $D'$ is either a bidirected path or an antidiracted path. If $D'$ is a bidirected path, then $H'$ is the disjoint union of two bidirected paths. If $D'$ is an antidirected path, then $H'$ is the disjoint union of two antidirected paths (both are copies of $D'$). Hence, for every vertex $x\in X_{v_1}$, there exists a partial transversal $T_x$ of $(X,H)$ with $\dom(T_x:D) = V(P)$,  $x\in T_x$, and $T_x$ is independent in $H$.
Hence $H[T_x]$ is strictly $f$-degenerate. Note that $f(y)\not=(0,0)$ for every vertex $y\in V(H)$ (by Claim~\ref{claim:no_0_0}). Let $X_{v_1}=\{x_1,x_2\}$. Then $X_{v_0}=\{y_1,y_2\}$ such that, in $D$, $y_1$ is adjacent to $x_1$ but not to $x_2$. For $i\in \{1,2\}$, let $\cK^i=(D^i,X^i,H^i,f^i)=\cK/T_{x_i}$. From \ref{enum:proof:thm2:A} and Proposition~\ref{prop_reduction_general} it then follows that $\cK^i$ is 
is degree-feasible and uncolorable. By \ref{enum:proof:thm2:B} and \ref{enum:proof:thm2:C}, this implies that $\cK^i$ is constructible. Note that $D^i=D-V(P)$ and $(X^i,H^i)=(X,H)/(D-V(P))$. Since $D-V(P)$ is 2-connected, this implies that $\cK^i$ is an \upshape{M}-, \upshape{K}-, \upshape{C}-, or  \upshape{A}-configuration. Since the vertex $u$ is, in $D$, not adjacent to any vertex in $V(P)$, we have $f^i(x_u)=f(x_u)$ for every vertex $x_u\in X_u$. Since $|X_u|=2$, we obtain that $\cK^i$ is not an M-configuration. 	

Assume that $\cK^1$ is a \upshape{K}- or a \upshape{C}-configuration. This implies that there are two transversals of $(X^1,H^1)$, say $T_1$ and $T_2$, such that $(X^1,H^1)/T_j$ is a saturated cover of $D^1 =D-V(P)$ and $f^1$ is constant on $T_j$ (for $j\in \{1,2\}$). Since $D^2=D^1$ and $(X^2,H^2)=(X^1,H^1)$, we obtain that $\cK^2$ is a \upshape{K}- or a \upshape{C}-configuration, too. Since $f^1(x_u)=f^2(x_u)=f(x_u)\not=(0,0)$ for every $x_u\in X_u$, we obtain that $f^2=f^1$. Recall that $y_1\in X_{v_0}$ is, in $D$, adjacent to $x_1$ but not to $x_2$. This implies that $f(y_1)\not= f^1(y_1)$ (by Claim~\ref{claim:reduction})
and $f^2(y_1)=f(y_1)$. Hence, we have $f(y_1)\not= f^1(y_1)=f^2(y_1)=f(y_1)$, which is impossible. 

Hence both  $\cK^1$ and $\cK^2$ are \upshape{A}-configurations. In particular, this implies that both $f^1$ and $f^2$ are valued in $\{(0,1),(1,0)\}$. Again, we have $f^2(y_1) = f(y_1)$ as $y_1$ is not adjacent to $x_2$, and $f^1(y_1) \neq f(y_1)$. Hence either $f^2(y_1)$ or $f^1(y_1)$ does not belong to $\{(0,1),(1,0)\}$, yielding the contradiction. This completes the proof. 
\end{proofclaim}

In what follows, let $S(D)=\set{v\in V(D)}{D-v \mbox{ is 2-connected}}$. By Claim~\ref{claim:Dremains}, we have $S(D)\not= \ems$. Furthermore, if $x\in V(H)$, we denote by $u(x)$ the unique vertex $u\in V(D)$ such that $x\in X_u$. If $v\in S(D)$ and $x\in X_v$, then let $\cK_x=\cK/(v,x)$ and let $f_x$ be the vertex function of $\cK_x$. Then we have $\cK_x=(D',X',H',f_x)$ where $D'=D-v$ and $(X',H')=(X,H)/D'$. By Claim~\ref{claim:reduction}, $\cK_x$ is a constructible configuration. Since $D'$ is $2$-connected, we have that $\cK_x$ is an \upshape{M}-, \upshape{K}-, \upshape{C}-, or  \upshape{A}-configuration. Let us first establish some basic facts about these configurations. 

\begin{claim}
\label{claim:Dreduced}
Let $v\in S(D)$ an arbitrary vertex and let $y\in X_v$ be an arbitrary vertex. Furthermore, let $D'=D-v$ and $(X',H')=(X,H)/D'$.
Then the following statements hold:
\begin{enumerate}[label={\rm (\alph*)}]
\item If $f_y(x)=(0,0)$ for a vertex $x\in V(H')$, then we have
\begin{enumerate}[label=\arabic*.]
      \item $f(x)=(a_H(x,y),a_H(y,x))=(a_D(u(x),v),a_D(v,u(x))$ and $x\in N_H^+(y) \cup N_H^-(y)$.
      \item $f_y(x')=f(x')\not= (0,0)$ for every vertex $x'\in X_{u(x)}\sm \{x\}$ and $f(X_{u(x)}\sm \{x\})=d_{D'}(u(x))$. 
    \end{enumerate}
\label{enum:claim:Dreduced:a}
\item \label{enum:claim:Dreduced:b}
If $\cK_y$ is an \upshape{M}-configuration, then $r=2$ and $V(H)$ is the disjoint union of two transversals $T$ and $T'$ of $(X,H)$ such that $T'=\suc(f)$, booth $(X',H')/T$ and $(X',H')/T'$ are saturated covers of $D'$ and 
$$f(x)=\begin{cases}
            d_{D'}(u(x)) & \text{if } x\in T, \\
			(a_D(u(x),v),a_D(v,u(x)) & \text{if } x\in T'.
		\end{cases}$$
As a consequence, we have $d_{D'}(u)=(a_D(u,v),a_D(v,u))$ for every vertex $u\in V(D')$.
\end{enumerate}
\end{claim}
\begin{proofclaim}
For the proof of \ref{enum:claim:Dreduced:a}, assume that $f_y(x)=(0,0)$, and let $u=u(x)$. Using Claim~\ref{claim:reduction}, we obtain that
$f(x)=f_y(x)+(a_H(x,y),a_H(y,x))=(a_H(x,y),a_H(y,x))$. By Claim~\ref{claim:perfect_matching}, this implies that $f(x)=(a_D(u,v),a_D(v,u))$. This proves the first statement. Furthermore, for every vertex $x'\in X_u\sm \{x\}$, we have $x'\not\in N_H^+(y) \cup N_H^-(y)$ and so $f(x')=f_y(x')+(a_H(x',z),a_H(z,x'))=f_y(x')$. Using \eqref{eq_degree-sum}, we then obtain that
$$f(X_u\sm \{x\})=f(X_u)-f(x)=d_D(u)-(a_D(u,v),a_D(v,u))=d_{D'}(u).$$
Thus \ref{enum:claim:Dreduced:a} is proved. For the proof of \ref{enum:claim:Dreduced:b} assume that $\cK_y=(D',X',H',f_y)$ is an \upshape{M}-configuration. Then there is a transversal $T$ of $(X',H')$ such that $(X',H')/T$ is a saturated cover of $D'$ and $f(x)=d_{D'}(u(x))$ for every $x\in T$ and $T'=V(H')\sm T=\suc(f)$. From \ref{enum:claim:Dreduced:a} it then follows that $T'$ is a transversal of $(X',H')$, $r=2$, and $f(x)=(a_D(u(x),v),a_D(v,u(x))$ if $x\in T'$. Since $(X',H')/T$ is a saturated cover of $D'$ and $r=2$, it follows from Claim~\ref{claim:perfect_matching}\ref{enum:claim:perfect_matching:a} that $(X',H')/T'$ is also a saturated cover of $D'$. This proves \ref{enum:claim:Dreduced:b}
\end{proofclaim}

To complete the proof of Theorem~\ref{theorem:main2}, let $S^+(D)$ denote the set of vertices $v\in S(D)$ for which
	\begin{enumerate}[label=(\roman*)]
		\item $d_G(v)$ is minimum, and
		\item $d^+_D(v) + d^-_D(v)$ is minimum subject to~(i).
	\end{enumerate}
Clearly, we have $S^+(D)\not=\ems$ and $S^+(D)\subseteq S(D)$. For $n=|D|$, we obviously have $n\geq 4$. 
 
\begin{claim}
\label{claim:D-v=complete}
Let $v\in S^+(D)$ be an arbitrary vertex, and let $y\in X_v$ be an arbitrary vertex. Then $\cK_y$ is a \upshape{K}-configuration but no \upshape{M}-configuration.
\end{claim}
\begin{proofclaim}
Suppose this is false. Then there is a vertex $y\in X_v$ such that $\cK_y=(D',X',H',f_y)$ is an \upshape{M}-, \upshape{C}-, or  \upshape{A}-configuration. Assume first, that $\cK_y$ is an \upshape{M}-configuration. By Claim~\ref{claim:Dreduced}\ref{enum:claim:Dreduced:b}, we then have that $d_{D'}(u)=(a_D(u,v),a_D(v,u))$ for every vertex $u\in V(D')$. Since $D'=D-v$ is $2$-connected. this implies that $d_{D'}(u)=(1,1)$ for every vertex $u\in V(D')$. Consequently, $D'$ is a directed cycle with $|D'|\geq 3$ and $N_D^+(v) = N_D^-(v) = V(D - v)$. Then $S(D)=V(D)$. 
Let $u$ be an arbitrary vertex of $D'$. Then $d_G(v) \geq 3 \geq d_G(u)$ and $d^+_D(v)+d^-_D(v) \geq 6 > d^+_D(u)+d^-_D(u)$, a contradiction to the choice of $v$.

Consequently, $\cK_y$ is a \upshape{C}-, or  \upshape{A}-configuration. Then $G'=G-v$ is the underlying graph of $D'$, which is a cycle with $|G'|\geq 4$. If $d_G(v)=|G|-1$, then $S(D)=V(D)$ and $d_G(v)\geq 4>d_G(u)$ for every vertex $x\in V(D')$, contradicting the choice of $v$. Hence there is a vertex $u'\in V(D')$ such that $vu'\not\in E(G)$. This implies that $d_G(u')=2$. Then it follows from Claim~\ref{claim:vertexG} that $r=2$ and so $(X,H)$ is a $2$-uniform cover of $D$. Consequently, we have $f_y(x)=(1,1)$ for all $x\in V(H')$ or $f_y(x)\in \{(1,0),(0.1)\}$ for all $x\in V(H')$. There is a vertex $z\in X_v\sm \{y\}$. Then we have that $\cK_z$ is a \upshape{C}-, or  \upshape{A}-configuration, too. Since $u'v\not\in E(G)$, we have $E_H(X_v,X_{u'}) \cup E_H(X_{u'},X_v)=\ems$ and so $f_y(x)=f(x)=f_z(x)$ for every vertex $x\in X_{u'}$ (by Claim~\ref{claim:reduction}). This implies that both configurations $\cK_y$ and $\cK_x$ are of the same type and $f_y(x)=f_z(x)$ for every vertex $x\in V(H')$. Clearly, $v$ has a neighbor $u$ in $G$, Then it follows from Claim~\ref{claim:perfect_matching} that there is a vertex $x\in X_u$ such that $x\in N_H^+(y)\sm N_H^+(z)$ or $x\in N_H^-(y)\sm N_H^-(z)$. In the former case we obtain that $f^-_y(x)=f^-(x)-1$ and $f^-_z(x)=f^-(x)$, and in the latter case we obtain that $f^+_y(x)=f^+(x)-1$ and $f^+_z(x)=f^+(x)$ (by Claim~\ref{claim:reduction}). Hence, in both cases we obtain that $f_y(x)\not=\f_z(x)$, a contradiction. This proves the claim. 
\end{proofclaim} 
	
\begin{claim}
\label{claim:D-v=comA}
Let $v\in S^+(D)$ be an arbitrary vertex, let $D'=D-v$ and $(X',H')=(X,H)/D'$. Then the following statements hold:
\begin{enumerate}[label={\rm (\alph*)}]
\item $D'=D^\pm(K_{n-1})$ and $V(H)$ is the disjoint union of $r$ transversals $T_1, T_2, \ldots, T_r$ of $(X',H')$ such that $H[T_i]=D^\pm(K_n)$ for $i\in [1,r]$. We call $(T_1,T_2, \ldots, T_r)$ a \DF{feasible sequence} of transversals of $(X',H')$. 
\label{enum:claim:D-v=comA:a}
\item \label{enum:claim:D-v=comA:b}
For every vertex $y\in X_v$, there are integers $n_i^y\in \nat$ for $i\in [1,r]$ such that $f_y(x)=(n_i^y,n_i^y)$ for every vertex $x\in T_i$. Furthermore, $n_1^y+n_2^y+\cdots +n_r^y=n-2$ and $n_i^y=0$ for at most one index $i\in [1,r]$. We call $(n_1^y, n_2^y, \ldots, n_r^y)$ a \DF{$y$-feasible sequence} of integers.
\item \label{enum:claim:D-v=comA:c} $D=D^\pm(K_n)$, that is, $N_D^+(v)=N_D^-(v)=V(D-v)$.  
\end{enumerate}
\end{claim}
\begin{proofclaim}
We proof \ref{enum:claim:D-v=comA:a} and \ref{enum:claim:D-v=comA:b} simultaneously. So let $y\in X_v$ be an arbitrary vertex. Then $\cK_y=(D',X',H',f_y)$ is a a \upshape{K}-configuration but no \upshape{M}-configuration (by Claim~\ref{claim:D-v=complete}). Since $(X,H')$ is an $r$-uniform cover of $D'$, this implies that $D'=D^\pm(K_{n-1})$, and there exist transversals $T_1,T_2, \ldots ,T_r$ of $(X',H')$ and $r$ integers $n_1,n_2, \dots, n_r$ such that:
		\begin{itemize}
			\item $n_1+n_2+\cdots +n_r = n-2$,
			\item  for every   $i\in [1,r]$ and $x\in T_i$ we have $f_y(x) = (n_i,n_i)$, and
			\item if $n_i \neq 0$ then $H'[T_i] = D^\pm(K_{n-1})$.
		\end{itemize}
By Claim~\ref{claim:Dreduced}\ref{enum:claim:Dreduced:a}, we obtain that $n_i=0$ for at most one index $i$. Using Claim\ref{claim:perfect_matching}\ref{enum:claim:perfect_matching:a}, we then obtain that also in this case $H[T_i]=D^\pm(K_{n-1})$. This proves
\ref{enum:claim:D-v=comA:a} and \ref{enum:claim:D-v=comA:b}. Note that if we change the vertex $y\in X_v$ to another vertex $z\in X_v$, then the feasible sequence of transversals of $(X',H')$ remains the same, only the feasible sequence of integers can change.

The proof of \ref{enum:claim:D-v=comA:c} is by contradiction. First assume that there is a vertex $u\in N_D^+(v)\sm N_D^+(v)$. Using Claims \ref{claim:perfect_matching} and \ref{claim:reduction}, there are vertices $y,z\in X_v$ and $x\in X_u$ such that $x\in N_H^+(y)\sm N_H^-(z)$. Then $x\in T_i$ for some $i\in [1,r]$ and  we obtain that $(n_i^y,n_i^y)=f_y(x)=f(x)-(0,1)$ and $(n_i^z,n_i^z)=f_z(x)=f(x)$. which is impossible. Similarly, we can argue, if $u\in N_D^-(v)\sm N_D^-(v)$. Hence we have $N_D^+(v)=N_D^-(v)$. Now assume that $u\in V(D')$ is not adjacent to $v$ in $D$. Clearly, there is a vertex $w\in N_D^+(v)$. Hence there is a vertices $y,z\in X_v$ and $x_w\in X_w$ such that $x_w\in N_H^+(y)\sm N_H^+(z)$. Then $x_w\in T_i$ for some $i\in [1,r]$ and there is a unique vertex $x_u\in X_u \cap T_i$. Then we obtain that $(n_i^y,n_i^y)=f_y(x_u)=f(x_u)=f_z(x_u)=(n_i^z,n_i^z)$ and $(n_i^y,n_i^y)=f_y(x_w)\not= f_z(x_w)=(n_i^z,n_i^z)$, which is impossible. 
This proves \ref{enum:claim:D-v=comA:c}.
 \end{proofclaim}
 
That $D=D^\pm(K_n)$ follows from Claim~\ref{claim:D-v=comA}\ref{enum:claim:D-v=comA:c}. This implies that $S^+(D)=S(D)=V(D)$. By Claim~\ref{claim:D-v=comA}\ref{enum:claim:D-v=comA:a}, this implies that for every vertex $v\in V(D)$, the cover $(X,H)/(D-v)$ has a feasible sequence of transversals. Since $n\geq 4$, this implies that $V(H)$ is the disjoint union of $r$ transversals of $(X,H)$, say $T_1, T_2, \ldots, T_r$, such that $H[T_i]=D^\pm(K_n)$. Let $x,x'\in T_i$ be two vertices. Then there is a vertex $v\in V(D-u(x)-u(x'))$ and a vertex $y\in T_j\cap X_v$ for an index $j\in  [1,r] \sm \{i\}$. Then $(T_1\sm X_v, T_2\sm X_v, \ldots, T_r\sm X_v)$ is a feasible sequence of transversals of $(X,H)/(D-v)$ and we have $f(x)=f_y(x)=(n_i^y,n_i^y)=f_y(x')=f(x')$. This shows that, for every $i\in [1,r]$ there is an integer $n_i\in \nat$, such that $f(x)=(n_i,n_i)$ for every $x\in T_i$. By Claim~\ref{claim:no_0_0}, we have $n_i\not=0$ for every $i\in [1,r]$. Using \eqref{eq_degree-sum} for a vertex $v\in V(D)$, we obtain that
$$(n-1,n-1)=d_D(v)=f(X_v)=n_1+n_2+\cdots +n_r.$$
Consequently, $\cK=(D,X,H,f)$ is a \upshape{K}-configuration. This contradiction to \ref{enum:proof:thm2:B} completes the proof of Theorem~\ref{theorem:main2}. \qedhere
\end{proof}

\section{Proofs of Theorems \ref{theorem:mainDP} and \ref{theorem:mainLI}}
\label{sectsub:proofs:mainDP_LI}

This section is devoted to the proof of Theorems~\ref{theorem:mainDP} and \ref{theorem:mainLI}. We recall them here for the reader's convenience.

\theoremmainDP*

\theoremmainLI*

\begin{proof}[Proof of Theorems~\ref{theorem:mainDP} and \ref{theorem:mainLI}]
Let $\cD$ be a reliable digraph property, let $D$ be a digraph, let $(X,H)$ be a $\cD$-critical cover of $D$, let $B$ be a block of the low vertex subdigraph $D[V(D,X,H,\cD)]$, and let $D'=D-V(B)$.
In the case of Theorem~\ref{theorem:mainLI}, we assume that the cover $(X,H)$ is associated with the considered list assignment. 

If $|B|=1$ then $B=K_1$ and we are done; so assume that $|B|\geq 2$. Since $(X,H)$ is a $\cD$-critical cover of $D$ and $D'$ is a proper induced subdigraph of $D$, there is a partial $\cD$-transversal $T$ of $(X,H)$ such that $\dom(T:D)=V(D')$ and $H[T]\in \cD$.

Let $(X^B,H^B)=(X,H)/B$ be the restricted cover of $B$. We now define a vertex function $f:V(H') \to \nato^2$ as follows. For every vertex $u\in V(B)$ and every color $x\in X_u$, we set
\[
	f^-(x) = \max \Big(0,d^-(\cD) - |N^-_H(x) \cap T| \Big) \text{~~~and~~~} f^+(x) = \max\Big(0,d^+(\cD) - |N^+_H(x) \cap T|\Big).
\]
Then $\cK=(B,X^B,H^b,f)$ is a configuration and we distinguish two cases.

\medskip

\case{1}{$\cK$ is colorable.} 

Then $(X^B,H^B)$ has a transversal $T_B$ such that $H_B[T_B]$ is strictly $f$-degenerate.
Note that $H^B[T_B]=H[T_B]$ and $T'\cup T_B$ is a transversal of $(X,H)$. Since $(X,H)$ has no $\cD$-transversal, we have $H[T' \cup T_B]\not\in \cD$. Hence there is a subset $T_1 \subseteq T' \cup T_B$ such that $H[T_1]\in \CR(\cD)$ (by Proposition~\ref{prop:smooth}\ref{enum:prop:smooth:c}). Since $H[T']\in \cD$ and $\cD$ is hereditary, we have $T_1\cap T_B\not=\ems$. Since $H[T_B]$ is strictly $f$-degenerate, this implies that $\tilde{H}= H[T_1\cap T_B]$ contains a vertex $x$ such that $d^+_{\tilde{H}}(x)<f^+(x)$ or $d^-_{\tilde{H}}(x)<f^-(x)$. By definition, $x\in X_u$ for some low vertex $u\in V(B)$.
	If $d^+_{\tilde{H}}(x)<f^+(x)$, then $f^+(x)=d^+(\cD)-|N^+_H(x) \cap T|$, which gives
	\[
	d^+_{H[T_1]}(x)\leq d^+_{\tilde{H}}(x)+|N^+_H(x) \cap T|<f^+(x)+|N^+_H(x) \cap T|= d^+(\cD). 
	\]
	Since $H[T_1]\in \CR(\cD)$,
	this contradicts Proposition~\ref{prop:smooth}\ref{enum:prop:smooth:e}. If $d^-_{\tilde{H}}(x)<f^-(x)$, we obtain that $d^+_{H[T_1]}(x)<d^-(\cD)$ with a similar argument, again a contradiction. This settles the first case.
	
\medskip

\case{2}{$\cK$ is uncolorable.}

Recall that every vertex $u$ of $B$ is a low vertex of $D$, that is, we have $d_D(u)=d(\cD)\cdot |X_u|$. For a vertex $u\in V(B)$, we obtain from Proposition~\ref{prop:cover} that
	\[
	\sum_{x\in X_u} |N^-_H(x) \cap T| \leq d^-_D(u)-d_B^-(u) \text{~~~and~~~} \sum_{x\in X_u} |N^+_H(x) \cap T| \leq d^+_D(u)-d_B^+(u),
	\]
	which, together with the fact that $u$ is a low-vertex, leads to
	\begin{eqnarray*}
		\sum_{x\in X_u}f^-(x)&\geq & \sum_{x\in X_u} (d^-(\cD)- |N^-_H(x) \cap T|)\\
		 &\geq&d^-(\cD)\cdot |X_u|-\sum_{x\in X_u} |N^-_H(x) \cap T| \geq d_B^-(u),
	\end{eqnarray*}
	and similarly $\sum_{x\in X_u}f^+(x) \geq d_B^+(u)$. This implies that $\cK=(B,X^B,H^B,f)$ is an uncolorable degree-feasible configuration. By Theorem~\ref{theorem:main2} it then follows that $\cK$ is a constructible configuration. Since $B$ is a block, this implies that $B$ is an 
\upshape{M}-, \upshape{K}-, \upshape{C}-, or  \upshape{A}-configuration.
	
	If the cover $(X,H)$ is associated with a list assignment of $D$, then $\cK$ cannot be a an even \upshape{C}-configuration nor an \upshape{A}-configuration. Hence we are done if $\cK$ is a \upshape{K}-, \upshape{C}-, or  \upshape{A}-configuration.
	
	It remains to consider the case that $\cK$ is a M-configuration. Then there is a transversal $T$ of $(X^B,H^B)$ such that $(X^B,H^B)/T$ is a saturated cover of $D$, and for $u\in V(B)$ and $x\in X_u$, we have $f(x)=d_B(u)$ if $x\in T$ and $f(x)=(0,0)$ otherwise. In particular, $H^B[T]=H[T]$ is a copy of $B$. For every vertex $u\in V(B)$ there is a unique vertex $x_u\in T\cap X_u$. By definition of $f$, we then have
	\[
	d_B(u)=f(x_u)\leq d(\cD).
	\]
	 Hence $\De(B)\leq d(\cD)$. If $B\in \cD$, then $B$ is a $\cD$-brick and we are done. It remains to consider the case that $B\not\in \cD$. From Proposition~\ref{prop:smooth}\ref{enum:prop:smooth:c}, it then follows that $B$ contains an induced subdigraph $B'\in \CR(\cD)$. Then $\de(B') \geq d(\cD)$. Since $B$ is a block and $\De(B)\leq d(\cD)$, this implies that $B=B'$ and $B$ is an Eulerian digraph such that $d_B^+(u)=d_B^-(u)= d^+(\cD)=d^-(\cD)$. This completes the proof. \qedhere
\end{proof}

\section{Concluding Remarks}
\label{sect:concluding}

\begin{figure}[htbp]
	\centering
	\begin{tikzpicture}[>=Stealth, every node/.style={rectangle, scale=0.8}]
		\node[draw, align=center] (thm1) at (0,0) {{\bf Brooks' Theorem}\\Brooks~\cite{Brooks41}};
		\node[draw, align=center] (thm2) at (-5.3,2) {{\bf Degenerate Partitioning} \\ Bollob\'as and Manvel~\cite{Bollobas79}, Borodin~\cite{Borodin1976}};
		\node[draw, align=center] (thm3) at (0,2) {{\bf Degree-choosability} \\ Borodin~\cite{BorodinThesis}, Erd\H{o}s et al.~\cite{ERT79}};
		\node[draw, align=center] (thm4) at (5,2) {{\bf Directed Brooks' Theorem} \\ Mohar~\cite{Mohar2010}};
		
		\node[draw, align=center] (thm5) at (-5.3,4) {{\bf Variable degeneracy} \\ Borodin et al.~\cite{BorodinKT00}};
		\node[draw, align=center] (thm6) at (0,4) {{\bf DP-degree-colorability} \\ Bernshteyn et al.~\cite{BernshteynKP2017}};
		\node[draw, align=center] (thm7) at (5,4) {{\bf Degree-dichoosability} \\ Harutyunyan and Mohar~\cite{HarutMo2011}};
		
		\node[draw, align=center] (thm8) at (-5.3,6) {{\bf DP-Variable degeneracy} \\ Lu et al.~\cite{Lu2022}, Kostochka et al.~\cite{KostochkaSS2023}};
		\node[draw, align=center] (thm9) at (0,6) {{\bf Directed Variable degeneracy} \\ Bang-Jensen et al.~\cite{BangJensenSS2020}};
		\node[draw, align=center] (thm10) at (5,6) {{\bf Bivariable degeneracy} \\ Gonçalves et al.~\cite{GoncalvesPR2024}};

		\node[draw, align=center] (thm11) at (0,8) {{\bf DP-Bivariable degeneracy} \\ Theorem~\ref{main_theorem}};
		
		\draw[->] (thm2.south) to (thm1.north west);
		\draw[->] (thm3.south) to (thm1.north);
		\draw[->] (thm4.south) to (thm1.north east);
		\draw[->] (thm5.south) to (thm2.north);
		\draw[->] (thm5.south east) to (thm3.north west);
		\draw[->] (thm6.south) to (thm3.north);
		\draw[->] (thm7.south west) to (thm3.north east);
		\draw[->] (thm7.south) to (thm4.north);
		\draw[->] (thm8.south) to (thm5.north);
		\draw[->] (thm8.south east) to (thm6.north west);
		\draw[->] (thm9.south west) to (thm5.north east);
		\draw[->] (thm9.south east) to (thm7.north west);
		\draw[->] (thm10) to (thm9);
		\draw[->] (thm11.south east) to (thm10.north west);
		\draw[->] (thm11.south west) to (thm8.north east);
	\end{tikzpicture}
	\caption{Generalisations of Brooks' Theorem. An arrow from Theorem~A to Theorem~B illustrates the fact that Theorem~B is a particular case of Theorem~A. The result from~\cite{Lu2022} and~\cite{KostochkaSS2023} is derived from Theorem~\ref{main_theorem} by restricting it to configurations $(D,X,H,f)$ in which $D$ and $H$ are bidirected and $f^+=f^-$. The result from~\cite{GoncalvesPR2024} is derived by restricting to configurations $(D,X,H,f)$ in which $H$ is a disjoint union of copies of $D$. For a detailed explanation on other results' implications, the reader is referred to~\cite{GoncalvesPR2024}.}
	\label{fig:generalisations_Brooks}
\end{figure}
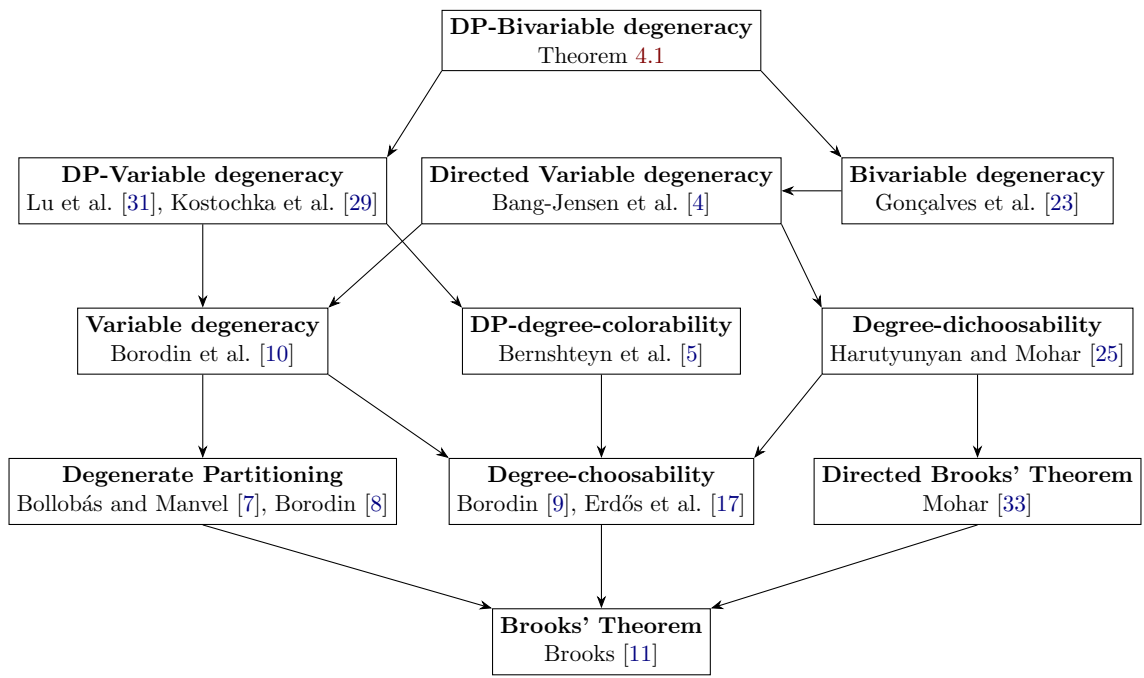

The partitioning and coloring of graphs under given degree constraints is a well-established area of graph theory that has attracted significant attention over the last 60 years. The concept of generalized graph coloring, in which vertices of the same color satisfy a given graph property, dates back to the late 1960s. Interested readers will find relevant information in \cite[Sections 2.5, 2.6 and 2.8]{StiebitzST2024}. Schweser \cite{Schweser2021} first pointed out that generalized (hyper)graph coloring is closely related to partition of (hyper)graphs into sub(hyper)graphs of prescribed variable degeneracy. For the ordinary coloring concept this was first observed by Borodin, Kostochka, and Toft \cite{BorodinKT00}. Generalized coloring of graphs with respect to the DP-chromatic number was first investigated by Kostochka, Schweser, and Stiebitz \cite{KostochkaSS2023}. 

\smallskip

With few exceptions, papers on digraph colorings deal with $\cAD$-colorings of digraphs, which are also known as \DF{acyclic colorings}. Acyclic colorings of digraphs were introduced in the late 1970s by Erd\H{os} and Neumann-Lara (see \cite{Erdos79} and \cite{NeumannLara82}). This is an appropriate extension of the ordinary coloring concept for graphs, as an acyclic coloring of a bidirected graph $D=D^\pm(G)$ induces a proper coloring of its underlying graph $G$, and vice versa. Hence, applying results about acyclic coloring of digraphs to bidirected graphs leads to result about proper coloring of graphs. As one might expect, many papers on acyclic coloring of digraphs seek to extend a coloring result for proper coloring of graphs to the class of digraphs. In this way, a large number of coloring results for graphs have found a counterpart for digraphs. It seems to be a difficult and often unsolvable problem to find counterparts of results about the partition of graphs satisfying degree constraints for the class of digraphs. Therefore, only a few results deal with the partition of digraphs that satisfy degree constraints. Most of these results are related to degeneracy (see e.g. \cite{BangJensenSS2020}, \cite{Golowich2016} and \cite{GoncalvesPR2024}). 

\smallskip

Our main result is Theorem~\ref{main_theorem}. If we apply this theorem to bidirected graphs and to symmetric vertex functions $f$ (i.e., $f^+=f^-$), we obtain the main result in \cite[Theorem 8]{KostochkaSS2023}. Note, however, that Theorem 8 is a result about multigraphs and not only about graphs. Extending Theorem~\ref{main_theorem} to multidigraphs would be worthwhile. On the other hand, Theorem~\ref{main_theorem} subsumes a large collection of different generalisations of Brooks' Theorem. Many of these results are summarized in Figure~\ref{fig:generalisations_Brooks}.

\smallskip

As shown in Section~\ref{sectsub:proofs}, Theorem~\ref{main_theorem} implies Theorems \ref{theorem:mainDP} and \ref{theorem:mainLI} immediately. 
Using these results, we easily obtain Corollaries~\ref{corollary:degree-choosable} and \ref{corollary:brooksDP}, and, moreover, Theorems \ref{theorem:brooksPD} and \ref{theorem:brooksLI}. Let us consider these results in case of $\cD=\cAD$. We denote by $T_2$ the digraph consisting of two vertices and one arc. Furthermore, let 
$$\cBR=\set{D^\pm(K_n)}{n\geq 1}
\cup \set{D^\pm(C_n)}{n\geq 3 \mbox{ odd}} \cup \{T_2\} \cup \set{\dic{n}}{n\geq 2}.$$ 
Using \eqref{equation:CR(AD)}, we obtain that $B$ is a $\cAD$-dibrick if and only if $B \in \cBR$. 

Then Theorem~\ref{theorem:brooksLI} gives the following result. If $D$ is a connected digraph, then $\lacn(D)\leq \max\{\De^+(D),\De^-(D)\}+1$ (by Corollary~\ref{corollary:brooksDP}\ref{enum:cor:brooksDP:c}) and equality holds if and only if $D=D^\pm(K_n)$ with $n\geq 1$, or $D=D^\pm(C_n)$ with $n\geq 3$ odd, or $D=\dic{n}$ with $n\geq 2$. This result was obtained by Harutyunyan and Mohar \cite{HarutMo2011}; the result for the dichromatic number was obtained by Mohar \cite{Mohar2010}. If we apply this result to bidirected graphs, we obtain Brooks' classical theorem for the chromatic number, respectively Brooks theorem for the list-chromatic number, which was first proved, independently, by Erd\H{o}s, Rubin, and Taylor~\cite{ERT79} and by Borodin~\cite{BorodinThesis}. 

Using Corollary~\ref{corollary:degree-choosable}\ref{enum:cor:degree-choosable:a}, we obtain the following result. Let $D$ be a connected digraph
with $|D|\geq 2$. If there exists a cover $(X,H)$ of $D$ such that $|X_v|\geq \max\{d_D^+(v),d_D^-(v)\}$ for every $v\in V(D)$ and $D$ has no $(\cAD,(X,H))$-coloring, then every block $B$ of $D$ satisfies $B=D^\pm(K_n)$ with $n\geq 2$, or $B=D^\pm(C_n)$ with $n\geq 3$, or $B=\dic{n}$ with $n\geq 2$. It is not difficult to show that the converse statement also holds. This result was obtained by Bang-Jensen \etal \cite[Theorem 15]{BangBSS2018}. The corresponding result for bidirected graphs was obtained by Bernshteyn, Kostochka, and Pron \cite{BernshteynKP2017}. 
If we apply Corollary~\ref{corollary:degree-choosable}\ref{enum:cor:degree-choosable:b} to bidirected graphs we obtain the following degree-choosable result. Let $G$ be a connected graph. Then there exists a list assignment $L$ of $G$ with $|L(v)|\geq d_G(v)$ for all $v\in V(G)$ such that $G$ has no $(\cO,L)$-coloring if and only if every block of $G$ is a complete graph or an odd cycle. This result was proved by  by Erd\H{o}s, Rubin, and Taylor~\cite{ERT79} and by Borodin~\cite{BorodinThesis}.

From Corollary~\ref{corollary:brooksDP} we obtain the following results. Let $D$ be a digraph such that 
$$D\in \Critaa{k+1} \cup \Critbb{k+1} \mbox{ (or, let } D\in \Critcc{k+1}\mbox{)}$$ 
with $k\geq 2$, and let $U=\set{v\in V(D)}{d_D^+(v)=d_D^-(v)=k}$. Then $\de(D)\geq (k,k)$ and every block $B$ of $D[U]$ satisfies $B\in \cBR$ (or, $B$ satisfies $B\in \cBR$, or $B=D^\pm(C_n)$ with $n\geq 4$ even, or $B$ is an antidirected cycle). For the dichromatic number, this result was first proved by Bang-Jensen \etal \cite[Theorem 15]{BangJensenBSS2019}. As pointed out in \cite{BangJensenBSS2019}, there are digraphs $D\in \Critaa{k+1}$ such that $T_2$ is a block of the low vertex subdigraph $D[U]$. If $D=D^\pm(G)$, then it follows from \eqref{equation:D=G} 
that $D\in \Critaa{k+1}$ if and only if $G\in \Crit(\cO,\cn,k+1)$. Hence we obtain Gallai's fundamental result from \cite{Gallai63a} about the block structure of $\cn$-critical graphs. If  $G\in \Crit(\cO,\cn,k+1)$ with $k\geq 3$ and $U=\set{v\in V(G)}{d_G(v)=k}$, then $\de(G)\geq k$ and every block of the low vertex subgraph $G[U]$ is a complete graph or an odd cycle. Gallai proved also that if $H$ is a connected graph with $\De(H)\leq k$ and each block of $H$ is a complete graph or an odd cycle, then there is a graph $G\in \Crit(\cO,\cn,k+1)$ such that the low vertex subgraph $G[U]=H$. Note that it follows from K\"onig's characterization of bipartite graphs that $\Crit(\cO,\cn,3)=\set{C_n}{n\geq 3 \mbox{ is odd}}.$

\smallskip

Coloring of digraphs under variable degeneracy constraints were first studied by Bang-Jensen, Schweser and Stiebitz \cite{BangJensenSS2020} and by 
Gon\c{c}alves, Picasarri-Arrieta, and Reinald \cite{GoncalvesPR2024}. They investigated the following coloring problem for digraphs. Let $p\in \nat$ be a fixed integer, and let $(D,\vf)$ be a pair such that $D$ is a digraph and $\vf=(f_1,f_2, \ldots, f_p)$ is a \DF{vector function} of $D$, i.e., $f_i:V(D) \to \nato^2$ for $i\in [1,p]$. We say that $(D,\vf)$ is \DF{colorable} if there is a coloring $\f$ of $D$ with color set $C=[1,p]$ such that $D[\fin(i)]$ is strictly $f_i$ degenerate for every color $i\in C$, for otherwise we say that $(D,\vf)$ is \DF{uncolorable}. We call $(D,\vf)$ \DF{degree-feasible} if every vertex $v\in V(D)$ satisfies
$$\sum_{i=1}^p f_i(v)\geq d_D(v).$$
A  characterization of uncolorable degree-feasible pairs $(D,\vf)$ whose underlying digraph $D$ is connected was obtained by Gon\c{c}alves, Picasarri-Arrieta, and Reinald \cite{GoncalvesPR2024}. The case when $\vf$ is symmetric (i.e. $f_i^+=f_i^-$) was handled in \cite{BangJensenSS2020}. 

For the class of digraphs a characterization of uncolorable degree-feasible pairs $(D,\vf)$ can be easily obtained from Theorem~\ref{main_theorem}. To this end, let $(D,\vf)$ be a pair such that $D$ is a connected digraph and $\vf$ is a vector function of $D$. 
Let $(X,H)=C(D,L)$ be the cover of $D$ associated to the constant list assignment $L: v \mapsto C=[1,p]$, that is, $X_v=\{v\}\times C$ for all $v\in V(D)$, and for two distinct vertices $(u,i)$ and $(v,j)$ of $H$ we have
$$a_H((u,i),(v,j))=
\left\{ \begin{array}{ll}
a_D(u,v) & \mbox{\rm if } i=j,\\
0 & \mbox{\rm if } i\not=j.
\end{array}
\right.$$
Note that $(X,H)$ is indeed a cover of $D$. Define a vertex function $g: V(H)\to \nato^2$  by $g(u,i)=f_i(u)$ for $u\in V(D)$ and $i\in C$. Then $\cK=(D,X,H,g)$ is a configuration. Furthermore, it is easy to check that $(D,\vf)$ is degree-feasible if and only if $\cK=(D,X,H,g)$ is degree feasible; and $(D,\vf)$ is colorable if and only if $\cK$ is colorable. Hence Theorem~\ref{main_theorem} yields a  characterization of uncolorable degree-feasible pairs $(D,\vf)$, provided that $D$ is a connected digraph. This is exactly the characterization given in \cite{GoncalvesPR2024}. If $(D,\vf)$ is an uncolorable degree-feasible pair and $D$ is a block, then it follows from Theorem~\ref{main_theorem} that $(D,\vf)$ satisfies one of the following three conditions:
\begin{itemize}
  \item $D$ is a block and there is a color $i\in C$ such that, for every vertex $v\in V(D)$, we have $f_i(v)=d_D(v)$ and $f_j(v)=(0,0)$ whenever $j\in C\sm \{i\}$.
  \item $D=D^\pm(K_n)$ for an integer $n\geq 1$, and the functions $f_1, f_2, \ldots, f_p$ are all constant and symmetric, and for every vertex $v\in V(D)$ we have $f_1(v)+f_2(v)+\cdots+ f_p(v)=(n-1,n-1)$.
  \item $D=D^\pm(C_n)$ for an odd integer $n\geq 5$, and the functions $f_1, f_2, \ldots, f_p$ are all constant and equal to $(0,0)$ except
  two functions that are constant and equal to $(1,1)$.
\end{itemize}
Note that since $(X,H)=C(D,L)$, the configuration $\cK=(D,X,H,g)$ can never be an even \upshape{C}-configuration nor an \upshape{A}-configuration. 

\smallskip

Given a configuration, it is easy to decide if it is degree-feasible. It also seems likely, that our proof of 
Theorem~\ref{main_theorem} can be transformed into a polynomial-time algorithm where the input is a degree-feasible configuration $\cK=(D,X,H,f)$ and the output is either a transversal $T$ of $(X,H)$ such that $H[T]$ is strictly $f$-degenerate or the blocks of $D$ showing that $\cK$ is a constructible configuration.

\end{document}